\documentclass[a4paper,fleqn]{cas-sc}

\usepackage[authoryear]{natbib}

\usepackage{algorithm}
\usepackage{algorithmicx}
\usepackage{algpseudocode}
\usepackage{eurosym}
\usepackage{placeins} 

\usepackage{color}
\usepackage[acronym]{glossaries}
\makenoidxglossaries

\def\tsc#1{\csdef{#1}{\textsc{\lowercase{#1}}\xspace}}
\tsc{WGM}
\tsc{QE}

\newcommand{\astar}{A$^*$}
\newcommand{\spoplus}{\acrshort{spo}+}

\newacronym{dfl}{DFL}{decision-focused learning}
\newacronym{evrp}{EVRP}{Electric Vehicle Routing Problem}
\newacronym{vrp}{VRP}{Vehicle Routing Problem}
\newacronym{lns}{LNS}{large neighborhood search}
\newacronym{milp}{MILP}{mixed-integer linear program}
\newacronym{pfl}{PFL}{prediction-focused learning}
\newacronym{spo}{SPO}{Smart Predict-then-Optimize}

\begin{document}
\let\WriteBookmarks\relax
\def\floatpagepagefraction{1}
\def\textpagefraction{.001}

\shorttitle{}    

\shortauthors{}  

\title [mode = title]{Carbon-Aware Electric Vehicle Routing Problem with Decision-Focused Learning}  

\author[1,2]{Jaume Reixach}[orcid=0009-0002-0305-9270]
\credit{Conceptualization, Methodology, Software, Validation, Formal analysis, Investigation, Data Curation, Writing -- original draft, Writing -- review \& editing}

\author[1,2]{Wanru Chen}[orcid=0000-0002-1864-2373]
\credit{Conceptualization, Methodology, Investigation, Writing -- review \& editing}

\author[1]{Jiaming Wu}[orcid=0000-0002-0235-4246]
\cormark[1]
\ead{jiaming.wu@chalmers.se}
\credit{Conceptualization, Methodology, Supervision, Writing -- review \& editing}

\author[2]{Balázs Kulcsár}[orcid=0000-0002-3688-1108]
\credit{Conceptualization, Methodology, Supervision, Writing -- review \& editing}

\affiliation[1]{organization={Architecture and Civil Engineering, Chalmers University of Technology},
            city={Gothenburg},
            country={Sweden}}

\affiliation[2]{organization={Electrical Engineering, Chalmers University of Technology},
            city={Gothenburg},
            country={Sweden}}

\cortext[1]{Corresponding author}

\begin{abstract}
Electric vehicles produce no tailpipe emissions, but charging them produces emissions that depend on the time-varying generation mix of the electricity consumed. We study an Electric Vehicle Routing Problem with a carbon-aware objective that accounts for charging-related emissions alongside conventional routing costs. The problem determines routing, charging, and waiting decisions, allowing charging demand to be shifted toward periods with lower grid carbon intensity. Because these decisions depend on forecasts of future carbon intensity, we propose a decision-focused learning approach so that the training of the model employed for predicting the carbon intensity of the grid over time is solver-aware. In particular, we train a linear predictor using approximate \spoplus{} directions and employ a time-limited Large Neighborhood Search followed by exact postprocessing as the underlying solver. Experiments on realistic Germany-based instances show that explicitly considering charging-related emissions reduces realized emissions by 28.76\% on average relative to carbon-blind optimization, although at the expense of higher non-carbon operating cost. Moreover, decision-focused training reduces realized total cost by 2.30\% on average relative to conventional prediction-focused learning and performs better on 75.5\% of the test instances, despite substantially lower pointwise predictive accuracy. Finally, we show that its average advantage remains positive across varied solver configurations and carbon-cost weights.
\end{abstract}

\begin{highlights}
    \item We formulate a carbon-aware EVRP with time-varying grid carbon intensity.
    \item We embed a time-limited LNS and exact postprocessing in approximate SPO+ training.
    \item Carbon-aware routing reduces charging-related emissions by 28.76\% on average.
    \item Decision-focused learning reduces realized cost by 2.30\% on average.
    \item Its benefits persist across solver configurations and carbon-cost weights.
\end{highlights}


\maketitle
\thispagestyle{plain}
\pagestyle{plain}
\glsresetall

\section{Introduction}\label{}

Road transport electrification is an important strategy for reducing greenhouse gas emissions.  In this context, the EU Electrification Action Plan \citep{europeancommission2026electrification} proposes increasing electricity’s share of the EU’s final energy consumption to 46\% by 2040. Transport is particularly relevant to this transition, as it accounts for almost one third of final energy consumption in the EU, approximately 75\% of which is attributable to road transport \citep{eurostat2026transport}. Electrifying road transport can therefore play an important role in achieving the EU’s broader electrification objectives. However, although electric vehicles produce no tailpipe emissions, their operation still generates indirect emissions through the production of electricity used for charging. Reducing these emissions requires attention to both electricity generation and how vehicles charge.

Charging-related emissions depend on both the amount of electricity consumed and its carbon intensity, which varies over time with the electricity generation mix. Consequently, two routing solutions with similar energy consumption can produce substantially different emissions. An energy-minimizing solution is therefore not necessarily emission-minimizing, since charging can be shifted toward periods of lower grid carbon intensity. Exploiting this opportunity requires coordinating routing, charging, and waiting decisions while satisfying operational constraints. We therefore study a carbon-aware Electric Vehicle Routing Problem (EVRP) that balances conventional operating costs against charging-related emissions under time-varying grid carbon intensity. Since future carbon intensity is unknown when these decisions are made, this coordination also depends on the forecasts used to guide optimization.

The classical approach to solving such problems is via two sequential steps: first training a predictive model, and then treating the prediction output as an input to an algorithm for solving the optimization problem. This approach, known as \gls{pfl} \citep{mandi2024decision}, treats prediction and optimization as separate tasks. However, high prediction accuracy does not necessarily translate into better decisions, as prediction errors differ in their impacts on the downstream optimization problem. In contrast, \gls{dfl} incorporates the downstream optimization problem into model training to improve decision quality \citep{mandi2024decision}. In this paper, we develop a \gls{dfl}-based method for the \gls{evrp}, training a linear predictor of grid carbon intensity over the routing horizon to maximize the use of clean electricity for charging. In particular, we train this predictor using approximate \spoplus{} directions, employing a \gls{lns} followed by exact fixed-route optimization of charging quantities and waiting times as the downstream solution method.

The main contributions of this paper are as follows:
\begin{itemize}
    \item We formulate a carbon-aware \gls{evrp} that jointly optimizes routing, charging, and waiting under time-varying grid carbon intensity.

    \item We develop a scalable solution procedure combining \gls{lns} with A$^*$ for charging-path construction and exact fixed-route postprocessing. The exact postprocessing step refines the charging quantities and waiting times of the routes produced by the \gls{lns}.

    \item We integrate this solution procedure into an approximate \spoplus{} training loop, allowing the carbon-intensity predictor to be trained
    directly with respect to its downstream routing performance.

    \item We evaluate the approach on realistic instances based on the
    German road network and electricity data, demonstrating that \gls{dfl} reduces average realized total cost relative to \gls{pfl} and other baselines.
\end{itemize}

The rest of the paper is organized as follows. Section~\ref{sec:literature_review} reviews the literature relevant to the study. Section~\ref{sec:problem_description} defines the considered \gls{evrp} and presents its \gls{milp} formulation. Section~\ref{sec:methodology} introduces the proposed \gls{dfl} framework and its downstream solver, including the \gls{lns} and the exact fixed-route postprocessing step. Section~\ref{sec:experimental_evaluation} evaluates the methodology on realistic Germany-based instances and compares it with conventional \gls{pfl}. Finally, Section~\ref{sec:conclusions} summarizes the main conclusions and outlines directions for future research.

\section{Literature Review}\label{sec:literature_review}

We first review papers where the problem tackled is relevant to the one considered in this work. Then, we present an overview of the most relevant existing methods for \gls{dfl}. Finally, we review the most relevant applications of \gls{dfl}, particularly the \spoplus{} loss, to \glspl{vrp} and compare them with the application presented in this paper.

\subsection{Carbon-Aware Routing and Charging}

The \acrfull{evrp} is a variant of the classical \gls{vrp} in which the vehicles employed are electric. This variant has become popular in recent years as electric vehicles have been gaining prevalence in real-world routing applications. The main advantages of using electric vehicles in routing applications include lower running costs and lower pollution. Lower costs are due to three main reasons: (i) electricity is usually cheaper than fuel, (ii) electric drivetrains are more efficient, and (iii) electric drivetrains generally require less maintenance than combustion drivetrains. On the other hand, electric vehicles can also produce less pollution as they do not burn fuel during driving, and the electricity used for charging the batteries can be generated from low-polluting sources.

An important detail is that, if the generation of electricity is done using carbon-intensive sources, such as coal, oil or natural gas, the use of electric vehicles can lead to higher lifecycle greenhouse-gas emissions than simply employing efficient hybrid or gas-powered vehicles. This is highlighted in \cite{singh2024ensuring}, where the authors show that in some regions of the United States, electric vehicles lead to higher lifecycle greenhouse emissions than some efficient hybrid counterparts due to a high grid carbon intensity. This motivates the main idea of this paper, which is to account for the variability in the carbon intensity of the grid, with the aim of charging during lower carbon intensity periods and ultimately using cleaner electricity for charging.

Most \gls{evrp} papers focus on minimizing the amount of vehicles employed, the total energy consumption, the total travelling time, or a combination of these values \citep{kucukoglu2021electric}. To our knowledge, there are no \gls{evrp} papers focusing on directly minimizing the carbon intensity of the electricity used for charging. In the wider field of routing problems, some similar ideas have been applied to end-to-end routing, for example in \cite{su2025optimizing}, the authors consider the time-varying carbon intensity of the grid and present an optimization model to schedule the charging of electric trucks in long-haul end-to-end routes across the United States. Another related paper is \cite{jarvis2022smart}, which concerns the scheduling of electric bus charging, but instead of explicitly considering the carbon intensity of the grid over the routing horizon, they predict the excess in clean energy originating from wind generation and aim to maximize the charging done during periods with a surplus of clean energy.

One of the \gls{evrp} papers conceptually closest to ours is \cite{miyabe2025low}, in which the authors consider the coordination of charging with solar energy production. In particular, they consider time-varying local surplus power from solar generation associated with each charging station, aiming to maximize the use of solar power for charging. Our paper instead focuses on treating the carbon intensity of the whole grid as uncertain and on forecasting it before making routing decisions. We consider a single forecasting model with the task of predicting the carbon intensity of the grid over time, while in the compared work they use an independent forecasting model for each station, which gathers features from a local area around it. Moreover, in the present paper we apply \gls{dfl} to make the training process aware of the downstream routing solver.

\subsection{Decision-Focused Learning}
Many real-world optimization problems follow a similar structure, in which the formulation of the problem depends on a vector of uncertain values that need to be predicted. Most commonly, these uncertain values appear in the objective function. Classical examples of this setting include routing applications in transport where travel times are not known in advance and portfolio allocation problems, in which the future value of assets is also unknown and needs to be predicted instead.

The structure of these types of problems is generally denoted as Predict-then-Optimize and the classical way of solving them consists of treating the prediction and optimization steps as independent. As mentioned in the introduction, this approach is usually denoted as \acrfull{pfl}. Under this approach, the prediction model is trained independently of the optimization problem, generally via standard supervised learning aiming to maximize prediction accuracy.

\Acrfull{dfl} is an alternative to \gls{pfl} that consists of making the training of the prediction model aware of the optimization problem. The main idea is to train the prediction model so that the performance of the downstream optimization algorithm is maximized, without explicitly optimizing its pointwise predictive accuracy.

The motivation behind this approach lies in the fact that small prediction errors can lead to big drops in downstream performance. The main challenge in \gls{dfl} is that training requires information about how the solution of the downstream optimization problem changes with the predicted parameters. This solution mapping is generally available only implicitly through an optimization solver and for combinatorial optimization problems it is often discontinuous and piecewise constant, so its derivatives are either undefined or equal to zero almost everywhere.

Following the taxonomy by \citet{mandi2024decision}, gradient-based \gls{dfl} methods can be organized into four families that apply to different classes of downstream optimization problems and obtain informative gradients in different ways. Broadly, these approaches seek to minimize decision regret, which is the increase in the true objective value caused by using the decision induced by the predicted coefficients rather than an optimal decision under the true coefficients.

First, we have approaches based on analytical differentiation that differentiate optimality conditions or iterations of a solution algorithm of continuous optimization problems \citep{amos2017optnet, kotary2023backpropagation, domke2012generic}. The second consists of methods that replace the original optimization problem with a smooth approximation, so that the resulting solution mapping can be differentiated \citep{wilder2019melding, ferber2020mipaal}. The third consists of perturbation-based methods that estimate gradients from the solutions obtained with the initial parameters perturbed \citep{domke2010implicit, berthet2020learning, niepert2021implicit}. Finally, surrogate-loss methods replace decision regret with a tractable loss designed to promote downstream solution quality \citep{elmachtoub2022smart, mulamba2021contrastive, mandi2022decision, shah2022decision}. These approaches differ in their assumptions about the optimization problem, their need for access to solver internals, and the number and type of optimization calls required during training.

In general, there is no universally preferable \gls{dfl} method: the appropriate choice depends on the downstream problem, the predictive model, the available access to the optimizer, and the computational cost of solving it.

Among these approaches, the surrogate loss \spoplus{} \citep{elmachtoub2022smart} is particularly compatible with the problem considered in this work. This is because the uncertain carbon-intensity vector enters the downstream \gls{evrp} linearly through the amount of energy charged in each time period, while the feasible region is independent of this prediction. Furthermore, \spoplus{} requires only solutions under modified objective coefficients and not differentiation through the internal operations of the solver, which makes it suitable for a problem involving discrete routing decisions. Although alternative perturbation-based \gls{dfl} methods could also be applied, \spoplus{} provides a direct, theoretically motivated surrogate for decision regret while requiring comparatively few calls to the computationally expensive downstream solver.

The original \spoplus{} construction assumes access to an exact optimization oracle. In our setting, repeatedly solving the full \gls{evrp} to optimality during training would be computationally impractical. We therefore use an approximate downstream solver consisting of a time-limited \gls{lns} followed by exact fixed-route postprocessing. The resulting update directions are approximate \spoplus{} subgradients and their effectiveness should therefore be assessed empirically as the theoretical guarantees established under exact optimization are not necessarily inherited. Our experimental evaluation shows that training with these approximate subgradients effectively improves downstream solution quality relative to conventional \gls{pfl}.

\subsection{Decision-Focused Learning for Vehicle Routing}

Applications of \gls{dfl} to vehicle routing remain relatively limited compared with applications to simpler combinatorial problems such as shortest paths. The work most closely related to our setting is \cite{chu2023data}, which considers a last-mile delivery problem formulated as a capacitated \gls{vrp} with uncertain travel times. The authors train a linear travel-time prediction model using \spoplus{} and employ a simulated annealing heuristic to obtain the routing solutions required during training. This combination of \spoplus{}, vehicle routing, and a heuristic optimization oracle is particularly close to the approach adopted in this paper. In both cases, the use of a heuristic makes it possible to perform the repeated routing computations required by decision-focused training without solving an exact formulation at every iteration. The main differences concern the predicted quantity and the downstream problem: \cite{chu2023data} predict arc travel times for a last-mile CVRP and use simulated annealing, whereas we predict time-dependent grid carbon intensities for an \gls{evrp} with charging decisions and use a time-limited \gls{lns} followed by exact fixed-route postprocessing as the downstream solver.

Another direct application of \gls{dfl} to vehicle routing is presented by \cite{mandi2021data}. Their objective is to learn implicit preferences from historical routes, represented through predicted arc probabilities, and incorporate these into a CVRP. Rather than using \spoplus{}, they treat an exact mathematical-programming solver as a black box and obtain surrogate gradients through additional solves with perturbed objective coefficients. This work demonstrates the applicability of \gls{dfl} to a CVRP, but differs from ours both in its learning objective and in its reliance on exact optimization. It also illustrates the importance of the amount of available training data, since the authors report that the decision-focused approach does not improve out-of-sample route quality in their limited-data experiments.

\Gls{dfl} has also been applied to other related routing settings. \cite{tian2023prescriptive} consider a maritime routing problem in which port-state-control officers must select and sequence ship inspections based on predictions of ship deficiencies. Since repeatedly solving the NP-hard routing problem during training is computationally expensive, the authors use surrogate losses based on noise-contrastive estimation and a pool of previously generated feasible solutions, thereby reducing the need for reoptimization. Although this is not a vehicle-delivery problem and does not use \spoplus{}, it highlights a central difficulty shared with our setting: incorporating a hard routing problem into the learning loop requires controlling the computational cost of repeated optimization.

Taken together, these studies show that \gls{dfl} can improve routing decisions, while also emphasizing the need for scalable downstream optimization.

\section{Problem Description}\label{sec:problem_description}

\begin{figure}[
    pos=h,
    width=0.7\textwidth,
    align=\centering
]
    \includegraphics[width=0.5\textwidth]
        {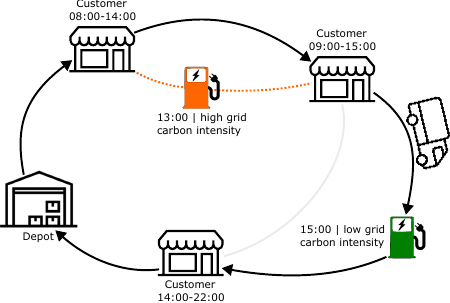}
    \caption{Example of a carbon-aware routing decision.}
    \label{fig:problem_illustration}
\end{figure}

We consider a carbon-aware electric vehicle routing problem with time windows. A homogeneous fleet serves customers from a central depot while adhering to capacity, battery, time window, and routing horizon constraints. Vehicles leave the depot fully charged and can visit charging stations, where they may wait before charging if beneficial. Upon returning to the depot, vehicles recharge to full capacity. Routing, charging, and waiting decisions are jointly optimized to minimize operational costs and indirect charging emissions. As grid carbon intensity varies over time, the timing of charging can affect the realized emissions even within solutions that consume similar amounts of energy. Figure~\ref{fig:problem_illustration} illustrates how time windows and the time-varying grid carbon intensity can affect charging-station selection. A tight customer time window prevents the vehicle from waiting at the first charging station. The vehicle therefore prefers a later detour over charging immediately at a higher grid carbon intensity.

We denote by \(0\) and \(n+1\) the depot at the beginning and end of a route, respectively, \(C\) the set of customers, and \(F\) the set of charging-station copies. Using copies allows a physical charging station to be visited multiple times while retaining a separate arrival time, battery state, and charging decision for each visit. The node set is \(V=\{0,n+1\}\cup C\cup F\). We assume a complete directed graph over \(V\), excluding self-loops, arcs entering \(0\), and arcs leaving \(n+1\). For every admissible arc \((i,j)\), \(\tau_{ij}\) and \(e_{ij}\) denote its travel time and energy consumption, respectively. Each customer \(i\in C\) has demand \(q_i\), service duration \(s_i\), and a time window associated \([T_i^L,T_i^U]\). Charging stations and depot nodes have no delivery demand, so we set \(q_i=0\) for \(i\in F\cup\{0,n+1\}\). All vehicles have delivery capacity \(Q_{\max}\), charging rate \(c_r\), maximum battery capacity \(y_{\max}\), and minimum permissible state of charge \(y_{\min}\).

The routing horizon \([0,t_{\max}]\) is partitioned into time buckets \(H\), where bucket \(h\in H\) corresponds to the interval \([\ell_h,u_h)\) and has grid carbon intensity \(CI_h\), measured in \(\mathrm{kgCO_2eq/kWh}\). The objective coefficients \(C_{\mathrm{veh}}\), \(C_{\mathrm{time}}\), and \(C_e\) denote the fixed cost per route, the cost of time, and the electricity price, respectively. Charging-related emissions are converted into monetary cost using the carbon-cost weight \(p_{\mathrm{CO_2}}\), measured in \(\text{€}/\mathrm{kgCO_2eq}\).

We formulate the problem as a \acrfull{milp}. The binary variable \(x_{ij}\) indicates whether arc \((i,j)\) is traversed. The continuous variables \(a_i\), \(b_i\), and \(d_i\) denote arrival, service or charging start, and departure times at node \(i\), respectively. The variables \(y_i\), \(L_i\), \(c_i\), and \(w_i\) represent the state of charge upon arrival, cumulative delivered load, energy charged, and waiting time. The binary variable \(z_{ih}\) indicates whether charging at station copy \(i\) starts in time bucket \(h\), while \(\phi_{ih}\) denotes the corresponding amount of energy priced using that bucket. For a route whose final arc is \((i,n+1)\), the variables \(b_i^D\), \(w_i^D\), and \(c_i^D\) describe the start time, waiting time, and amount of terminal recharging at the depot, while \(z_{ih}^D\) and \(\phi_{ih}^D\) assign this recharge to a time bucket. The constants $M_t$ and $M_e$ denote sufficiently large values; in practice, they can be set to  $M_t = t_{\max} + \max_{i,j\in V} \tau_{i,j}$ for the timing constraints and to $M_e = y_{\max} + \max_{i,j\in V}e_{i,j}$ for the battery constraints, or tightened further using arc-specific travel-time and energy bounds. To simplify the notation, we use the convention that arc variables and any summations involving them include only admissible ordered pairs. 

The objective minimizes a weighted combination of vehicle, travel-time, waiting, charging-time, electricity, and charging-related emission costs, as shown in Equation~\eqref{eq:cevrp_objective}. The first term counts the number of nonempty routes as the number of arcs leaving the start depot. The second term charges travel time over all selected arcs. The third term penalizes waiting at customers, charging stations, and the terminal depot, as well as charging time at intermediate stations. The final two terms price the energy charged at intermediate charging stations and at the terminal depot, using both the electricity price $C_e$ and the carbon-cost weight $p_{\mathrm{CO2}}$ multiplied by the carbon intensity $CI_h$ of the selected time bucket.
\begin{align}
\min \quad
& C_{\mathrm{veh}}\sum_{j\in V} x_{0j}
+ C_{\mathrm{time}}
  \sum_{i,j\in V}
  \tau_{ij}x_{ij}
\nonumber\\
& + C_{\mathrm{time}}
    \left(
    \sum_{i\in C\cup F} w_i
    + \sum_{i\in C\cup F} w_i^D
    + \sum_{i\in F}\frac{c_i}{c_r}
    \right)
\nonumber\\
& + \sum_{i\in F}\sum_{h\in H}
    \phi_{ih}\left(C_e+p_{\mathrm{CO2}}CI_h\right)
  + \sum_{i\in C\cup F}\sum_{h\in H}
    \phi^D_{ih}\left(C_e+p_{\mathrm{CO2}}CI_h\right).
\label{eq:cevrp_objective}
\end{align}

Note that a charging event is assigned to a single time bucket. Accordingly, the bucket duration is assumed to be at least as long as the maximum charging time. This discretization, aimed at simplifying the model, nevertheless remains approximate, since a charging event starting near the end of a bucket would in reality extend into the subsequent bucket. However, the approximation affects only the portion of the charge extending into that subsequent bucket and is expected to introduce limited error because grid carbon intensity varies continuously over time and therefore generally changes only moderately between contiguous time buckets.

Route connectivity and customer coverage are enforced by the following constraints. Each customer must have exactly one outgoing and one incoming selected arc, ensuring that every customer is visited exactly once. Charging stations are optional: if a copy is entered, it must also be left, and each copy can be used at most once. The final constraint balances the number of routes leaving the start depot with the number of routes arriving at the end depot.
\begin{align}
& \sum_{j\in V}x_{ij}=1
    && \forall i\in C
\label{eq:visit_customer_out}\\
& \sum_{j\in V}x_{ji}=1
    && \forall i\in C
\label{eq:visit_customer_in}\\
& \sum_{j\in V}x_{ij}
  =
  \sum_{j\in V}x_{ji}
    && \forall i\in F
\label{eq:station_flow}\\
& \sum_{j\in V}x_{ij}\le 1
    && \forall i\in F
\label{eq:station_once}\\
& \sum_{j\in V}x_{0j}
  =
  \sum_{i\in V}x_{i,n+1}.
\label{eq:depot_balance}
\end{align}

The next constraints track the cumulative demand delivered along each route. The cumulative delivered load is zero at the start depot. At a customer, it must be at least the demand of that customer and cannot exceed the vehicle's capacity. At a charging station, it is bounded by the same vehicle capacity but its lower bound is 0 as $q_i = 0$ if $i\in F$. The propagation constraint increases the cumulative delivered load by $q_j$ whenever arc $(i,j)$ is selected; again, since $q_j=0$ for charging stations, the same expression also applies when the next visited node is a charging station.
\begin{align}
& L_0=0
\label{eq:load_start}\\
& q_i\le L_i\le Q_{\max}
    && \forall i\in C\cup F
\label{eq:load_customer_bounds}\\
& L_j \ge L_i + q_j - Q_{\max}(1-x_{ij})
    && \forall i\in V,\ j\in C\cup F.
\label{eq:load_propagation} \\
& L_j \leq L_i + q_j + Q_{\max}(1-x_{ij})
    && \forall i\in V,\ j\in C\cup F.
\label{eq:load_propagation_2}
\end{align}

The timing constraints initialize the route at the start depot. When an arc is selected, the arrival time at the next customer or charging station must equal the departure time plus travel time; this is enforced with two big-$M$ inequalities. Service or charging may only start after arrival. Customer service must start inside the corresponding time window. Departure times equal service completion times at customers and charging completion times at charging stations. Waiting is defined as the difference between service or charging start time and arrival time. Finally, any route returning from node $i$ to the end depot must do so within the routing horizon $t_{\max}$.
\begin{align}
& a_0=d_0=0
\label{eq:time_start}\\
& a_j \ge d_i+\tau_{ij}-M_t(1-x_{ij})
    && \forall i\in V,\ j\in C\cup F
\label{eq:arrival_lower}\\
& a_j \le d_i+\tau_{ij}+M_t(1-x_{ij})
    && \forall i\in V,\ j\in C\cup F
\label{eq:arrival_upper}\\
& b_i\ge a_i
    && \forall i\in C\cup F
\label{eq:start_after_arrival}\\
& T_i^L\le b_i\le T_i^U
    && \forall i\in C
\label{eq:customer_time_windows}\\
& d_i=b_i+s_i
    && \forall i\in C
\label{eq:customer_departure}\\
& d_i=b_i+\frac{c_i}{c_r}
    && \forall i\in F
\label{eq:station_departure}\\
& w_i = b_i-a_i
    && \forall i\in C\cup F
\label{eq:waiting_definition}\\
& d_i+\tau_{i,n+1}\le t_{\max}+M_t(1-x_{i,n+1})
    && \forall i\in C\cup F.
\label{eq:return_horizon}
\end{align}

The battery constraints initialize each route with a full battery at the start depot. For selected arcs ending at customers or charging stations, the arrival state of charge is the previous state of charge plus any charge received at the origin, minus the energy required to traverse the arc. The reserve constraint ensures that every selected arc can be traversed while maintaining at least $y_{\min}$ energy. Battery states are bounded by the minimum reserve and the maximum battery capacity. Non-charging nodes cannot receive charge, whereas charging stations may charge up to the remaining battery capacity.
\begin{align}
& y_0=y_{\max}
\label{eq:battery_start}\\
& y_j \le y_i+c_i-e_{ij}+M_e(1-x_{ij})
    && \forall i\in V,\ j\in C\cup F
\label{eq:battery_upper}\\
& y_j \ge y_i+c_i-e_{ij}-M_e(1-x_{ij})
    && \forall i\in V,\ j\in C\cup F
\label{eq:battery_lower}\\
& y_i+c_i-e_{ij}\ge y_{\min}-M_e(1-x_{ij})
    && \forall i,j\in V
\label{eq:battery_reserve}\\
& y_{\min}\le y_i\le y_{\max}
    && \forall i\in C\cup F
\label{eq:battery_bounds}\\
& c_i=0
    && \forall i\in C\cup\{0\}
\label{eq:no_customer_charge}\\
& 0\le c_i\le y_{\max}-y_i
    && \forall i\in F.
\label{eq:station_charge_bounds}
\end{align}

Charging at intermediate stations is linked to the time-varying carbon intensity through the following constraints, which assign every visited charging station to exactly one time bucket. The lower and upper bucket bounds force the charging start time $b_i$ to lie in the selected interval $[\ell_h,u_h)$. Here, $\varepsilon>0$ is a small enough constant that ensures the time intervals do not overlap. The variable $\phi_{ih}$ can be positive only for the selected bucket, and the bucket-level charged energy equals the total amount charged at the charging station.
\begin{align}
& \sum_{h\in H}z_{ih}
  =
  \sum_{j\in V}x_{ij}
    && \forall i\in F
\label{eq:station_bucket_select}\\
& \ell_h z_{ih}\le b_i
    && \forall i\in F,\ h\in H
\label{eq:station_bucket_lower}\\
& b_i\le u_h-\varepsilon+M_t(1-z_{ih})
    && \forall i\in F,\ h\in H
\label{eq:station_bucket_upper}\\
& 0\le \phi_{ih}\le y_{\max}z_{ih}
    && \forall i\in F,\ h\in H
\label{eq:station_phi_bounds}\\
& c_i=\sum_{h\in H}\phi_{ih}
    && \forall i\in F.
\label{eq:station_charge_accounting}
\end{align}

Terminal-depot waiting and recharging are represented separately as follows. When a route returns from node $i$ to the end depot, the model computes the amount of energy required to recharge the vehicle to full after arrival. The first two inequalities define this quantity if arc $(i,n+1)$ is selected, and the third constraint forces it to zero otherwise. The next four constraints define the terminal recharge-start time $b_i^D$ as the depot arrival time $d_i+\tau_{i,n+1}$ plus the optional waiting time $w_i^D$. Both variables are forced to zero when the corresponding return arc is not selected. The remaining constraints assign the terminal recharge to the time bucket containing $b_i^D$. Thus, a vehicle may wait after returning to the depot and begin recharging in a later time bucket, with this waiting penalized in the objective.
\begin{align}
& c_i^D \ge y_{\max}-y_i-c_i+e_{i,n+1}-M_e(1-x_{i,n+1})
    && \forall i\in C\cup F
\label{eq:depot_charge_lower}\\
& c_i^D \le y_{\max}-y_i-c_i+e_{i,n+1}+M_e(1-x_{i,n+1})
    && \forall i\in C\cup F
\label{eq:depot_charge_upper}\\
& 0\le c_i^D\le y_{\max}x_{i,n+1}
    && \forall i\in C\cup F
\label{eq:depot_charge_bounds}\\
& 0\le b_i^D\le t_{\max}x_{i,n+1}
    && \forall i\in C\cup F
\label{eq:depot_start_bounds}\\
& 0\le w_i^D\le t_{\max}x_{i,n+1}
    && \forall i\in C\cup F
\label{eq:depot_wait_bounds}\\
& b_i^D \ge d_i+\tau_{i,n+1}+w_i^D-M_t(1-x_{i,n+1})
    && \forall i\in C\cup F
\label{eq:depot_start_lower}\\
& b_i^D \le d_i+\tau_{i,n+1}+w_i^D+M_t(1-x_{i,n+1})
    && \forall i\in C\cup F
\label{eq:depot_start_upper}\\
& \sum_{h\in H}z^D_{ih}=x_{i,n+1}
    && \forall i\in C\cup F
\label{eq:depot_bucket_select}\\
& \ell_h z^D_{ih}\le b_i^D
    && \forall i\in C\cup F,\ h\in H
\label{eq:depot_bucket_lower}\\
& b_i^D\le u_h-\varepsilon+M_t(1-z^D_{ih})
    && \forall i\in C\cup F,\ h\in H
\label{eq:depot_bucket_upper}\\
& 0\le \phi^D_{ih}\le y_{\max}z^D_{ih}
    && \forall i\in C\cup F,\ h\in H
\label{eq:depot_phi_bounds}\\
& c_i^D=\sum_{h\in H}\phi^D_{ih}
    && \forall i\in C\cup F.
\label{eq:depot_charge_accounting}
\end{align}

Finally, the domains of the decision variables are defined below.
\begin{align}
& x_{ij}\in\{0,1\}
    && \forall i,j\in V
\label{eq:x_domain}\\
& z_{ih}\in\{0,1\}
    && \forall i\in F,\ h\in H
\label{eq:z_domain}\\
& z^D_{ih}\in\{0,1\}
    && \forall i\in C\cup F,\ h\in H
\label{eq:z_depot_domain}\\
& 0 \leq a_i,b_i,d_i,w_i \leq t_{\max}
    && \forall i\in V
\label{eq:timing_domains}\\
& y_i,L_i\ge 0
    && \forall i\in V
\label{eq:continuous_node_domains}\\
& b_i^D,w_i^D,c_i,c_i^D\ge 0
    && \forall i\in C\cup F
\label{eq:continuous_charge_domains}\\
& \phi_{ih}\ge 0
    && \forall i\in F,\ h\in H.
\label{eq:continuous_phi_domains1} \\
& \phi_{ih}^D\ge 0
    && \forall i\in C\cup F,\ h\in H.
\label{eq:continuous_phi_domains2}
\end{align}

\section{Methodology}\label{sec:methodology}

\begin{figure}[pos=htbp]
    \centering
    \includegraphics[width=\textwidth]{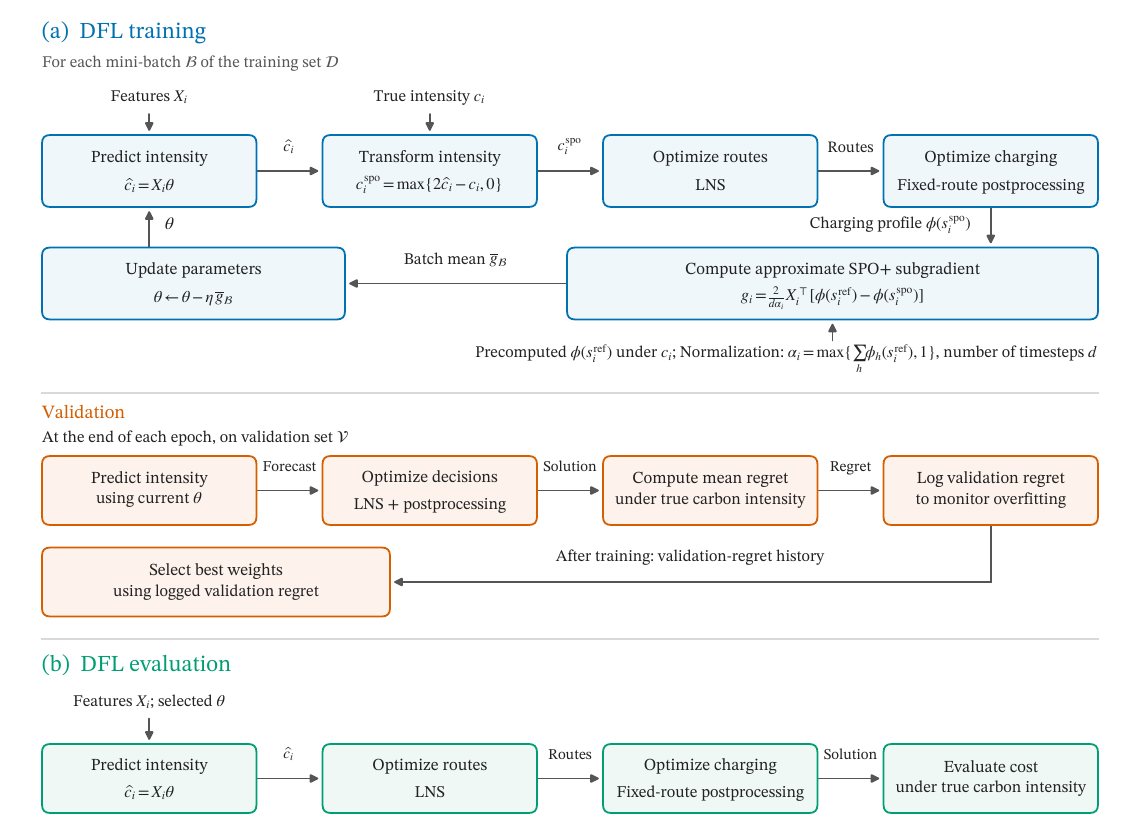}
    \caption{Overview of the proposed methodology: (a) decision-focused
training, validation, and model selection; (b) evaluation using the
trained predictor.}
    \label{fig:methodology_overview}
\end{figure}

This section presents the proposed decision-focused learning framework and its downstream optimization algorithm. Figure~\ref{fig:methodology_overview} summarizes the complete methodology. A linear model predicts grid carbon intensity over the routing horizon, and these predictions parameterize the carbon-related terms of the \gls{evrp} objective. A time-limited \gls{lns} followed by exact fixed-route postprocessing is used as the downstream solver. As already mentioned, this approximate downstream solver is necessary because repeatedly solving the full \gls{milp} to optimality for every training instance and parameter update would be computationally impractical for the realistic problem instances we are interested in.

During training, this solver is applied under transformed carbon-intensity vectors to construct approximate \spoplus{} subgradients. At the end of each epoch, the current predictor is evaluated on the validation set by computing a validation regret. This value is the average relative difference between the costs of the solutions obtained from the current predictions and those of the corresponding precomputed reference solutions, with both evaluated under the true carbon intensities. Validation regret is used for model selection.

During evaluation, the trained predictor supplies forecasts to the same algorithm, and the resulting solutions are evaluated under the realized carbon intensities. 

\subsection{Approximate \spoplus{} Training}
\label{subsection:training_loop}

As mentioned in the introduction, the \spoplus{} loss falls into the family of \gls{dfl} methods that differentiate a surrogate loss rather than the optimization mapping itself. In particular, the \spoplus{} loss is a surrogate of the so-called true \gls{spo} loss \citep{elmachtoub2022smart}. Consider an optimization problem of the form

\begin{equation}
\begin{aligned}
z^*(c) := \min_{w} \quad & c^T w \\
\text{s.t.} \quad & w \in S,
\end{aligned}
\label{eq:optimization_problem}
\end{equation}
where $w\in \mathbb{R}^d$ are the decision variables, $S\subseteq \mathbb{R}^d$ is the feasible region, and $c\in\mathbb{R}^d$ is the objective-coefficient vector, which corresponds to the carbon intensity forecast in our setting. We denote by $W^*(c) = \arg \min_{w\in S}\{c^T w\}$ the set of optimal solutions to (\ref{eq:optimization_problem}) for a cost vector $c$. The \gls{spo} loss, which the \spoplus{} loss is a surrogate of, is the following one:

\begin{equation}
    \ell_\text{SPO}(\hat{c}, c) = \max_{w\in W^*(\hat{c})}\{c^T w\} - z^*(c)
    \label{eq:spo_loss}
\end{equation}

The goal is to find a prediction $\hat{c}$ so that the \gls{spo} loss (\ref{eq:spo_loss}) is minimized, meaning that the true cost of the optimal solutions under $\hat{c}$ is as close as possible to the true cost of the optimal solutions under $c$. Note that the maximum $w$ over $W^*(\hat{c})$ is taken, as there might be multiple optimal solutions for $\hat{c}$ that have different values under $c$. The issue with this loss function is that it is hard to minimize directly, as it may not even be continuous in $\hat{c}$. To tackle this, the paper authors derive the \spoplus{} loss, which is a tractable upper bound for it:

\begin{equation}
    \ell_{\text{SPO}+}(\hat{c}, c) = \max_{w\in S}\{(c - 2\hat{c})^T w\} + 2\hat{c}^T w^*(c) - z^*(c).
    \label{eq:spo_loss_general}
\end{equation}

In this expression, $w^*(\cdot)$ denotes the solution returned by an optimization oracle. Although the maximum over $w\in S$ is generally non-differentiable, one can easily compute subgradients of (\ref{eq:spo_loss_general}), which allows us to employ mini-batch subgradient descent, as we do in this paper. 

In particular, we now assume that a linear predictor $\hat{c} = X\theta$ is used, and that we have a dataset consisting of the features and true values of $c$ for $n$ problem instances: $\{(X_i, c_i)\}_{i = 1}^n$. Following the Empirical Risk Minimization (ERM) principle, the goal is to minimize the mean loss over the dataset. That is, to solve

\[\min_{\theta\in \mathbb{R}^{p\times 1}} \frac{1}{n} \sum_{i = 1}^n \ell_{SPO+}(X_i\theta, c_i),\]

where $\theta\in\mathbb{R}^{p\times 1}$ are the parameters of the linear regression, and \(X_i\in\mathbb{R}^{d\times p}\) and \(c_i\in\mathbb{R}^d\) are the features and the true carbon intensity over the time steps of the \(i\)-th problem instance, respectively. Accordingly, \(d\) is the number of time steps and \(p\) the number of features associated with each time step. In this setting and following the authors, one can see that a valid subgradient of $ \ell_{SPO+}(X_i\theta, c_i)$ is

\begin{equation}
2X_i^T\big(w^*(c_i) - w^*(2X_i\theta - c_i)\big).
\label{eq:subgradient}
\end{equation}

The theoretical appeal of \spoplus{} is that, for a fixed realization $c$, it is a convex function of the prediction $\hat{c}$ and an upper bound on the \gls{spo} loss. Moreover, under sufficient conditions and in the population and unrestricted-prediction setting, \cite{elmachtoub2022smart} establish Fisher consistency of the \spoplus{} loss with respect to the \gls{spo} loss, meaning minimizing the expected \spoplus{} loss also minimizes the expected \gls{spo} loss.

An important detail is that these properties assume that \(w^*\) solves the problem exactly. Our use of an approximate method produces an approximate subgradient instead and therefore the theoretical guarantees established for exact optimization do not directly carry over. Below, we make the resulting subgradient of our setting explicit.

Let \(s\in\mathcal{F}\) denote a feasible \gls{evrp} solution and let \(\phi_h(s)\) be the total amount of energy charged during time bucket \(h\), including charging at intermediate stations and terminal depot recharging. After collecting all objective terms that do not depend on carbon intensity in \(G(s)\), the objective of a fixed solution can be written as

\[
G(s)+p_{\mathrm{CO2}}\sum_{h\in H}c_h\phi_h(s),
\]
where \(c\) is the carbon-intensity vector and \(p_{\mathrm{CO2}}\) is the carbon-cost weight. Notably, the objective is affine in the uncertain parameter \(c\), and the feasible set \(\mathcal{F}\) does not depend on its predicted value, placing the problem within the standard \spoplus{} setting.

For a prediction \(\widehat{c}=X\theta\), the \spoplus{} subgradient of Equation (\ref{eq:subgradient}) requires solving the downstream problem twice, once under the true cost vector \(c\) and once under the transformed cost vector \(2\widehat{c}-c\). Let \(s^{\mathrm{ref}}\) be the precomputed solution under the true carbon intensities and let \(s^{\mathrm{spo}}\) be the solution returned under the transformed vector. In the implementation, we define
\[
\alpha=\max\left\{\sum_{h\in H}\phi_h(s^{\mathrm{ref}}),1\right\}
\]
and use the normalized direction
\[
\widetilde g_\theta
=
\frac{2}{d\alpha}X^\top
\left(
\phi(s^{\mathrm{ref}})
-
\phi(s^{\mathrm{spo}})
\right).
\]
Dividing by \(\alpha\) prevents instances with longer routes or greater energy consumption from dominating an update solely because of their scale. Division by the number of time steps \(d\) keeps the direction on a comparable scale across routing horizons of different lengths. The subgradient also contains the constant carbon-cost weight \(p_{\mathrm{CO2}}\). Since the same carbon-cost weight is used for every training instance, this common positive factor only rescales all parameter updates and is therefore absorbed into the learning rate. The carbon-cost weight is nevertheless included in every \gls{lns} execution and thus still affects the routing and charging solutions used to construct the direction.

Reference solutions under the true carbon intensities can be precomputed, so in practice training requires one \gls{lns} and postprocessing execution per instance and gradient step. During the training of our experimental evaluation we precompute these solutions using the same time-limited \gls{lns} and postprocessing step, although with a higher time limit than the one used to compute \(s^{\mathrm{spo}}\).

Because the \gls{lns} is a time-limited heuristic rather than an exact optimization algorithm, \(s^{\mathrm{ref}}\) and \(s^{\mathrm{spo}}\) are approximate solutions. The resulting direction should therefore be interpreted as an approximate \spoplus{} subgradient associated with the time-limited \gls{lns} and the postprocessing. Furthermore, in the implementation we set negative components of the transformed carbon-intensity vector to zero because the embedded \astar{} search assumes nonnegative costs. This can be seen in line 6 of Algorithm~\ref{alg:Training_loop}. This projection provides a solver-compatible approximation but adds another layer of approximation as it departs from the exact \spoplus{} transformation whenever it is active. As reported in the experimental evaluation, this clipping affected less than 0.01\% of all transformed vector components in our setting.

Algorithm \ref{alg:Training_loop} presents the training loop, which applies mini-batch subgradient descent using the approximate \spoplus{} subgradient above. \(\mathcal{D}\) and \(\mathcal{V}\) are the training and validation sets, respectively, which contain one tuple per problem instance. Each tuple contains (i) features for every time step, (ii) the true carbon intensity vector over the time steps of the problem instance, and (iii) a precomputed solution under the true carbon intensity vector.

\begin{algorithm}[t]
\caption{Approximate \spoplus{} training loop}
\begin{algorithmic}[1]
    \Statex \textbf{Input}: Training set
        $\mathcal{D}=\{(X_i,\mathbf{c}_i,s_i^{\mathrm{ref}})\}_{i=1}^n$,
        validation set $\mathcal{V}$, initial parameters $\theta$, learning
        rate $\eta$, number of epochs $E$, batch size $m$
    \Statex \textbf{Output}: Epoch history
    $\{(\theta_e,R_e)\}_{e=1}^{E}$
    \For{$e=1,\ldots,E$}
        \State Randomly shuffle $\mathcal{D}$
        \For{each mini-batch $\mathcal{B}\subseteq\mathcal{D}$ of size at most $m$}
            \For{each instance $i\in\mathcal{B}$}
                \State $\widehat{\mathbf{c}}_i\gets X_i\theta$
                \State $\mathbf{c}_i^{\mathrm{spo}}\gets
                    \max\{2\widehat{\mathbf{c}}_i-\mathbf{c}_i,\mathbf{0}\}$
                    componentwise
                \State $s_i^{\mathrm{spo}}\gets
                    \textsf{apply\_postprocessing}(\textsf{LNS}(\mathbf{c}_i^{\mathrm{spo}}))$
                \State $\alpha_i\gets
                    \max\{\sum_{h\in H}\phi_h(s_i^{\mathrm{ref}}),1\}$
                \State $g_i\gets
                    \dfrac{2}{d\alpha_i}X_i^\top
                    \left(\boldsymbol{\phi}(s_i^{\mathrm{ref}})
                    -\boldsymbol{\phi}(s_i^{\mathrm{spo}})\right)$
            \EndFor
            \State $g\gets |\mathcal{B}|^{-1}
                \sum_{i\in\mathcal{B}}g_i$
            \State $\theta\gets\theta-\eta g$
        \EndFor
        \State Compute the average regret $\widehat R_{\mathcal V}(\theta)$ on $\mathcal{V}$
        \State Log $R_e\gets\widehat R_{\mathcal V}(\theta), \ \theta_e\gets\theta$
    \EndFor
    \State \textbf{return} $\{(\theta_e,R_e)\}_{e=1}^{E}$
\end{algorithmic}
\label{alg:Training_loop}
\end{algorithm}

For each epoch, the training set is shuffled and the parameters \(\theta\) are updated once per mini-batch. For an instance \(i\) in the mini-batch the transformed cost vector \(c_i^\mathrm{spo}\) is computed, with all of its negative components set to 0. Then, the \gls{lns} and postprocessing are executed with this transformed cost vector and the obtained solution is used together with the precomputed one \(s_i^\mathrm{ref}\) to compute the normalized approximate subgradient.

The validation set  \(\mathcal{V}\) is employed at the end of every epoch for computing a validation regret, which can be used to detect overfitting and select the trained weights. More specifically, for a validation instance \(i\), let \(s_i^\theta\) denote the solution obtained using the current prediction \(X_i\theta\), and let \(s_i^{\mathrm{ref}}\) denote the precomputed reference solution obtained using the true carbon-intensity vector \(c_i\). Writing \(\operatorname{cost}(s;c_i)\) for the total cost of solution \(s\) evaluated under \(c_i\), the percentage validation regret is:

\[
\widehat R_{\mathcal V}(\theta)
=
\frac{100}{|\mathcal V|}
\sum_{i\in\mathcal V}
\frac{\operatorname{cost}(s_i^\theta;c_i)-\operatorname{cost}(s_i^{\mathrm{ref}};c_i)}
{\operatorname{cost}(s_i^{\mathrm{ref}};c_i)}.
\]

Ideally, each reference solution \( s_i^\mathrm{ref}\) would be optimal for the true carbon-intensity vector, so that validation regret would be measured relative to a proven optimum. In our experimental evaluation, due to the computational impracticality of obtaining such optimal solutions, we precomputed them with the \gls{lns} and the postprocessing step under the same extended computational limit used for computing the reference solutions \(s^\text{ref}\) involved in the subgradient computation. Consequently, the reported validation regret should be interpreted as a measure of relative heuristic performance rather than exact optimization regret. In our implementation, we select the weights yielding the lowest five-epoch trailing mean validation regret. That is, we select \(\theta_{e^*}\), where \(e^*=\arg\min_{e\in\{5,\ldots,E\}}\frac{1}{5}\sum_{j=e-4}^{e}R_j\). The trailing mean is used to reduce the effect of the variability introduced by the stochastic solver.

\subsection{Large Neighborhood Search}
\Acrfull{lns} is a metaheuristic that has been widely applied to vehicle routing and other transportation problems \citep{pisinger2018large}. The main idea of the algorithm is to probabilistically destroy and repair a solution at every iteration with the goal of improving its quality over time. After application of the destroy and repair operators, the solution obtained is given a certain probability of replacing the solution kept by the algorithm, which will be used in the next iteration. Thus, solutions of a higher quality are given a higher replacement probability. 

The \gls{lns} operators destroy and repair parts of a solution by removing and reinserting a subset of the solution components forming it. In \glspl{vrp}, this usually corresponds to removing and reinserting customers, which is exactly our setting. Our \gls{lns} employs two of the most established repair and destroy operators, a regret-$k$ repair and a Shaw destroy \citep{pisinger2018large}.

Algorithm \ref{alg:LNS} provides pseudocode for the complete \gls{lns}. After each destroy-repair operation, the algorithm decides whether to accept the obtained solution as incumbent using simulated annealing \citep{kirkpatrick1983optimization}.

With more detail, if $\Delta_{\text{cost}} = S_{\text{cost}}^{\text{cand}} - S_\text{cost}$ is the difference between the cost of the new solution $S^\text{cand}$ and the solution $S$ currently being kept by the \gls{lns}, simulated annealing gives the following probability of accepting $S^\text{cand}$, and hence setting $S = S^\text{cand}$:

\begin{equation}
p = 
\begin{cases} 
1, & \text{if } \Delta_\text{cost} \leq 0, \\
e^{-\Delta_\text{cost} / T},  & \text{if } \Delta_\text{cost} > 0. 
\end{cases}
\end{equation}

Here, $T$ is the so-called temperature, which controls how likely a worse solution is to be accepted. This value is initialized at $T_\text{init}$ and decreased at each iteration by multiplying the so-called cooling factor $c_f$. This way, a worse solution is less likely to be accepted as the iterations progress, with the goal of escaping local minima early on. Notably, both $T_\text{init} > 0$ and $c_f \in (0,1)$ are algorithm parameters.

\begin{algorithm}[t]
\caption{Large Neighborhood Search}
\begin{algorithmic}[1]
    \State $ S \gets \textsf{generate\_initial\_solution}()$
    \State $ S^\text{bsf} \gets S$
    \State $T \gets T_\text{init}$
    \While{time limit not reached} 
      \State $S^\text{cand} \gets S$
      \State $\textsf{apply\_destroy}(S^\text{cand})$
      \State $\textsf{apply\_repair}(S^\text{cand})$
      \State $\Delta_{\text{cost}} \gets S^\text{cand}_{\text{cost}} - S_{\text{cost}}$
      \State $\text{acceptance} \gets \textsf{simulated\_annealing\_decision}(\Delta_{\text{cost}}, T)$
      \If{\text{acceptance}}
        \State $S \gets S^\text{cand}$
        \If{$S_{\text{cost}} < S^\text{bsf}_{\text{cost}} $}
          \State $S^\text{bsf} \gets S$
        \EndIf
      \EndIf
      \State $T \gets T \cdot c_f$     
    \EndWhile
    \State{$\textbf{return } S^\text{bsf}$}
\end{algorithmic} 
\label{alg:LNS}
\end{algorithm}

\subsubsection{\astar{} charging-path construction}

When removing or adding customers within a route, the \gls{lns} needs to decide whether to insert charging-station stops within the new connections to ensure battery-feasible paths. To do this, we implement an \astar{} search \citep{hart1968formal} to construct a path between two consecutive route nodes. In the \astar{}, a search state is characterized by the current node, arrival time, and state of charge. From each state, the algorithm considers edges either to the target node or to charging stations reachable from the current node. An expansion is feasible only if the vehicle has sufficient energy to traverse the arc while maintaining the minimum battery level. The cost of a travel leg $(i,j)$ includes the driving-time cost $C_{\mathrm{time}}\tau_{ij}$. If the successor node is a charging station, the vehicle is charged to full in the path-construction phase; the added cost is the electricity cost, the carbon cost evaluated at the arrival-time bucket, and the charging-time cost. If the successor is the depot, the terminal recharge is also priced using the carbon intensity at the depot arrival time, but no terminal charging-time cost is added. Thus, for a charging amount $\Delta y$ starting in bucket $h$, the cost contribution is

\[
\Delta y\left(C_e+p_{\mathrm{CO2}}CI_h\right)
+
C_{\mathrm{time}}\frac{\Delta y}{c_r},
\]

with the last term omitted for terminal depot recharging. States that would violate the target-node time window upper bound are discarded. The priority of a state is $f=g+\eta$, where $g$ is the accumulated cost and the heuristic is

\[
\eta(i)=C_{\mathrm{time}}\frac{\mathrm{dist}_{\mathrm{geo}}(i,t)}{v_{\max}},
\]

where \(\mathrm{dist}_{\mathrm{geo}}(i, t)\) is the straight-line distance from the current node \(i\) to the target node \(t\) and \(v_{\max}\) is an upper bound on the vehicle speed. Thus, the heuristic is admissible as it is a lower bound on the travelling cost.

To limit the number of labels explored, the implementation retains only the state with the lowest accumulated cost at each node. Newly generated states whose cost is not lower than the current node label are discarded. Because arrival time and state of charge are not included in this comparison, this is a cost-based state-pruning heuristic rather than a complete multi-resource dominance rule. It reduces the computational cost of the repeated searches performed within the \gls{lns}, at the expense of potentially discarding a more expensive state with greater downstream time or battery flexibility.

\subsubsection{Construction and Neighborhood Operators}

The initial solution is generated through sequential insertion. The following process is iterated until all customers are visited: (i) a route visiting the farthest unvisited customer is built, and (ii) customers are iteratively inserted into it, by inserting, at each step, the customer leading to the lowest cost increase. The \astar{} is used for deciding the paths between nodes, so each route remains feasible with respect to battery, time-window, and capacity constraints.

At every \gls{lns} iteration, the Shaw destroy operator removes \(q\) customers that are related in terms of distance, arrival time, and demand. First, a random customer is added to the set $B$ of customers to remove. Then, a customer is added to $B$ iteratively, until $|B| = q$. This is done by first selecting a random customer $r$ from $B$ to serve as the reference customer, and then selecting a customer not in $B$ by giving a higher probability to the ones with a lower relatedness score with respect to $r$.

The relatedness score is measured by:
\[
R(r,c)=\phi_1 d(r,c)+\phi_2 t(r,c)+\phi_3 \ell(r,c),
\]
where each component is normalized by dividing it by the maximum value. Hereby, $d(r,c)$, $t(r,c)$, $\ell(r,c)$ are the normalized differences in geographical distance, arrival time, and demand, respectively. \(\phi_1\), \(\phi_2\) and \(\phi_3\) are real parameters that weight each relatedness component. Customers are ranked by increasing
relatedness score to the reference customer, and the selected rank is
\[
\rho=\left\lfloor |P|U^p\right\rfloor + 1,
\]
where \(|P|\) is the size of the current candidate pool, \(U\sim\mathcal{U}(0,1)\), and
\(p \geq 1\) is a real algorithm parameter that controls the degree of determinism. A value of \(p=1\) corresponds to random selection while larger values bias the selection more strongly towards customers related to the reference customer.

After selecting the removed customers, the affected routes are rebuilt by deleting the selected customer nodes and recomputing the consecutive subpaths between the remaining nodes using the \astar{}. Since these subpaths are optimized sequentially, a newly computed subpath may leave insufficient battery or time-window slack at its destination for the subsequent subpath. In this case, the procedure backtracks and repeats the preceding \astar search with a higher minimum battery level at arrival, up to a fixed retry limit. If failure occurs on the first subpath or the fixed retry limit is reached, the destroy operation is rejected and the \gls{lns} proceeds to the next iteration with the current solution unchanged. Otherwise, empty routes are discarded and the resulting feasible partial solution, together with the removed customers, is passed to the repair operator.

Removed customers are inserted using regret-\(k\) repair. For every unvisited customer \(c\), the algorithm identifies its feasible insertion alternatives and computes

\[r_c = \sum_{\ell = 2}^{\min (k, |\mathcal{I}_c|)} (\Delta_{\ell c} - \Delta_{1c}),\]

where \(\mathcal{I}_c\) is its set of feasible insertions and \(\Delta_{l c}\) the cost of its \(l\)-best insertion. The customer with the largest regret is inserted at its cheapest feasible position, and the process is iterated until all customers have been reinserted. Feasible insertions are evaluated using \astar{}, and creating a new singleton route is always considered. This option ensures reinsertion whenever the customer is individually feasible.

\subsection{Fixed-route Charging and Waiting Optimization}
\label{subsection:postprocessing}
An important observation is that the \gls{lns} as presented always assumes charging to full and does not allow deliberate waiting. To control the amount charged and the waiting time without introducing these decisions into every destroy and repair evaluation, we introduce a postprocessing step that we apply to the final solution produced by the \gls{lns}. This is important, as waiting, when compatible with the time-window constraints, can move a subsequent charge to a period with lower carbon intensity, while partial charging can shift energy consumption from high-carbon-intensity time buckets to lower-carbon-intensity ones.

Postprocessing fixes each route produced by the \gls{lns}, including its customer order, travelled arcs, and charging-station visits, and optimizes only its schedule and charging amounts. Our implementation solves one \gls{milp} independently for each route. It associates a waiting variable with every position of the fixed route: at customers this includes waiting required by the lower time-window bound and may also delay later charging, while at the terminal depot it may postpone the terminal recharge.

Consider a fixed route written as a sequence $r=(v_0,\ldots,v_m)$, where \(v_0\) and \(v_m\) represent the depot, and the sequence includes every intermediate visit to a charging station selected by the \astar{}. Let $P=\{0,\ldots,m\}$ be the set of route positions, $P^-=\{0,\ldots,m-1\}$ the positions with an outgoing leg, $P_C=\{p\in P:v_p\in C\}$ the customer positions, and $P_F=\{p\in P:v_p\in F\}$ the charging-station positions. For each $p\in P^-$, we define $\tau_p=\tau_{v_pv_{p+1}}$ and $e_p=e_{v_pv_{p+1}}$. As in the full formulation, each time bucket $h\in H$ has lower and upper bounds $\ell_h$ and $u_h$. 

The continuous variables $a_p$, $b_p$, and $d_p$ are the arrival,
activity start (service for customers, charging for the rest), and departure times at route position $p$; $w_p$ is the
waiting time; $y_p$ is the state of charge upon arrival; and $c_p$ is
the amount charged. For $p\in P_F$, the binary variable $z_{ph}$ selects
the time bucket in which charging starts, and $\phi_{ph}$ is the amount
charged and priced in that bucket. The corresponding variables for the
terminal depot recharge are denoted by $z_h^D$ and $\phi_h^D$. The
terminal recharge amount itself is represented by $c_m$.

The objective aims at minimizing the same weighted costs as the full formulation. Since the route, the charging-station visits, and all travelled legs are fixed, the vehicle and driving-time costs are constant and can be omitted from the optimization objective. The remaining terms minimize waiting time, charging time at intermediate stations, and the electricity and carbon costs of intermediate and terminal charging:
\begin{align}
\min \quad
& C_{\mathrm{time}}
  \left(
    \sum_{p\in P}w_p
    +\sum_{p\in P_F}\frac{c_p}{c_r}
  \right)
\nonumber\\
&+\sum_{p\in P_F}\sum_{h\in H}
  \phi_{ph}\left(C_e+p_{\mathrm{CO2}}CI_h\right)
 +\sum_{h\in H}
  \phi_h^D\left(C_e+p_{\mathrm{CO2}}CI_h\right).
\label{eq:post_objective}
\end{align}
As with the full \gls{milp}, the terminal recharge is priced but its duration is not included, since
the route is complete when the vehicle returns to the depot.

The following constraints determine the schedule of travel, service, charging, and waiting along the fixed route. The route starts at time zero. At every position, the activity may start
no earlier than arrival, and waiting is exactly the difference between
these times. Because the route is fixed, arrival at the next position is
given directly by the departure time at the current position plus the
known travel time. Customer service starts within its time window and
lasts for the prescribed service time. At charging stations, departure
occurs after the selected amount has been charged at charging speed $c_r$. At the
terminal depot, departure is identified with the start of terminal
recharging because its duration is outside the route horizon.
\begin{align}
& a_0=b_0=d_0=w_0=0
\label{eq:post_time_start}\\
& w_p=b_p-a_p
    && \forall p\in P
\label{eq:post_wait_definition}\\
& a_{p+1}=d_p+\tau_p
    && \forall p\in P^-
\label{eq:post_arrival_flow}\\
& T_{v_p}^L\le b_p\le T_{v_p}^U
    && \forall p\in P_C
\label{eq:post_customer_windows}\\
& d_p=b_p+s_{v_p},
  \quad c_p=0
    && \forall p\in P_C
\label{eq:post_customer_activity}\\
& d_p=b_p+\frac{c_p}{c_r}
    && \forall p\in P_F
\label{eq:post_station_departure}\\
& 0\le c_p\le y_{\max}-y_p
    && \forall p\in P_F
\label{eq:post_station_charge_bounds}\\
& d_m=b_m,
  \quad c_m=y_{\max}-y_m.
\label{eq:post_terminal_activity}
\end{align}

Battery charge is propagated along the fixed route as follows. The vehicle leaves the initial depot with a full battery. Along every leg, the arrival state of charge at the next position equals the
arrival state at the current position plus any energy charged there,
minus the energy required by the leg. Bounds on every $y_p$ enforce both
the minimum reserve and the battery capacity. Together with
Constraint~\eqref{eq:post_station_charge_bounds}, these constraints also
prevent charging beyond full capacity.
\begin{align}
& y_0=y_{\max},
  \quad c_0=0
\label{eq:post_battery_start}\\
& y_{p+1}=y_p+c_p-e_p
    && \forall p\in P^-
\label{eq:post_battery_flow}\\
& y_{\min}\le y_p\le y_{\max}
    && \forall p\in P.
\label{eq:post_battery_bounds}
\end{align}

Intermediate charging is assigned to carbon-intensity time buckets through the following constraints. As in the full \gls{milp}, each charging event is assigned entirely to the time bucket in which charging begins. The big-$M$ bounds place its charging start time in the selected bucket, while
$\phi_{ph}$ can be positive only for that bucket. The sum of these
bucket-level quantities equals the amount charged at the visit. As in the full \gls{milp}, $\varepsilon > 0$ is a small constant that prevents the overlap of time buckets.
\begin{align}
& \sum_{h\in H}z_{ph}=1
    && \forall p\in P_F
\label{eq:post_station_bucket_select}\\
& \ell_h z_{ph}\le b_p
    && \forall p\in P_F,\ h\in H
\label{eq:post_station_bucket_lower}\\
& b_p\le u_h - \varepsilon + M(1-z_{ph})
    && \forall p\in P_F,\ h\in H
\label{eq:post_station_bucket_upper}\\
& 0\le\phi_{ph}\le y_{\max}z_{ph}
    && \forall p\in P_F,\ h\in H
\label{eq:post_station_phi_bounds}\\
& c_p=\sum_{h\in H}\phi_{ph}
    && \forall p\in P_F.
\label{eq:post_station_charge_accounting}
\end{align}

The terminal recharge is assigned to the bucket in which it starts. The
same linking constraints used for intermediate stations ensure that
$\phi_h^D$ equals $c_m$ in the selected bucket and is zero in every other bucket. In our implementation, we set $M=t_{\max}$.

\begin{align}
& \sum_{h\in H}z_h^D=1
\label{eq:post_depot_bucket_select}\\
& \ell_h z_{h}^D\le b_m
    && \forall h\in H
\label{eq:post_depot_bucket_lower}\\
& b_m\le u_h - \varepsilon + M(1-z_h^D)
    && \forall h\in H
\label{eq:post_depot_bucket_upper}\\
& 0\le\phi_h^D\le y_{\max}z_h^D
    && \forall h\in H
\label{eq:post_depot_phi_bounds}\\
& c_m=\sum_{h\in H}\phi_h^D.
\label{eq:post_depot_charge_accounting}
\end{align}

Finally, the domains of the remaining variables are specified as follows. All time variables are bounded by the route horizon $t_{\max}$, and all
charging variables are bounded by the battery capacity $y_{\max}$. The only integer
variables are the time-bucket selectors:
\begin{align}
& z_{ph}\in\{0,1\}
    && \forall p\in P_F,\ h\in H
\label{eq:post_station_z_domain}\\
& z_h^D\in\{0,1\}
    && \forall h\in H
\label{eq:post_depot_z_domain}\\
& 0\le a_p,b_p,d_p,w_p\le t_{\max}
    && \forall p\in P
\label{eq:post_time_domains}
\end{align}

The complete solution is postprocessed by solving this model separately
for every nonempty route. After solving, the implementation reconstructs
the route schedule and recomputes its complete cost, including the fixed
vehicle and travel costs omitted from \eqref{eq:post_objective}.

\FloatBarrier

\section{Experimental Evaluation}
\label{sec:experimental_evaluation}

We evaluate the proposed methodology on synthetic \gls{evrp} instances constructed using real-world German road-network, charging-infrastructure, weather, and electricity-system data. Germany provides a relevant setting because its grid carbon intensity varies substantially over time as the generation mix changes. This temporal variation creates opportunities to reduce charging-related emissions by optimizing routing, charging, and waiting decisions.

We focus on long-haul routes spanning multiple days. The need for en-route charging, together with carbon-intensity variations across multiple days, makes this setting suitable for examining how charging locations, quantities, and timing can be coordinated to balance operating costs and emissions.

\subsection{Problem Instances}

The considered problem instances have a 72-hour routing horizon and 50 to 100 customers, distributed throughout Germany. Each customer is associated with a standardized demand between 5 and 30, a service time between 0.5 and 2 hours, and a time window type, with all three assignments done uniformly at random. We consider three types of time windows: (i) a flexible time window covering the entire routing horizon, \([0,72]\); (ii) a day-specific time window covering one of the intervals \([0,24]\), \([24,48]\), or \([48,72]\); and (iii) a working-hours time window within a selected day. For the latter, the opening and closing times are sampled uniformly from hours 5--9 and 17--21 of that day, respectively.

In order to make the time windows more realistic, we bias customers farther from the depot toward later days when the time windows are day-specific or fixed to the working hours of a day. To do this, we consider the straight-line distance between a customer and the depot, normalized by dividing by the maximum distance over all customers. If this value is smaller than 0.7 we assign probabilities 0.4, 0.3, and 0.3 to select the first, second, and third routing days, respectively. Otherwise, we associate the probabilities 0.1, 0.4, and 0.5, respectively.

The locations of depots, customers, and charging stations are based on OpenStreetMap (OSM) data \citep{openstreetmap}. Depot and customers are sampled uniformly from OSM elements with the relevant tags. Candidate depot locations are identified using the tags \texttt{industrial=depot} and \texttt{landuse=depot}, while candidate customer locations are identified using \texttt{shop=wholesale}, \texttt{industrial=warehouse}, and \texttt{building=warehouse}.

Candidate charging stations are identified using \texttt{amenity=charging\_station}, as well as \texttt{amenity=fuel} together with \texttt{charging\_station} or \texttt{fuel:electricity} set to \texttt{true}. Afterwards, two filtering steps are applied. First, we retain only stations located within 300 m of a major road, identified by the tags \texttt{highway=motorway}, \texttt{highway=motorway\_link}, \texttt{highway=trunk}, \texttt{highway=trunk\_link}, and \texttt{highway=primary}. This focuses the charging network on stations relevant to long-haul routing. Second, nearby stations are aggregated using DBSCAN \citep{ester1996density} with a 5~km neighborhood radius measured by routing distance. This reduces the concentration of stations around urban areas and merges multiple stations representing approximately the same charging location. Using routing distance instead of straight-line distance also avoids merging stations that are geographically close but located on opposite sides of a highway. These steps produce a set of 2,231 charging stations.

A further instance-specific filtering step is applied to reduce computational requirements. For each instance, we first retain only stations reachable from at least one customer or the depot. Germany is then divided into \(50\times50\)~km grid cells and, for each customer and the depot, only the closest charging station within each grid cell is retained. After this filtering, the generated instances contain an average of 2,041 charging stations.

\begin{figure}[pos=htbp]
    \centering
    \includegraphics[width=0.6\textwidth]{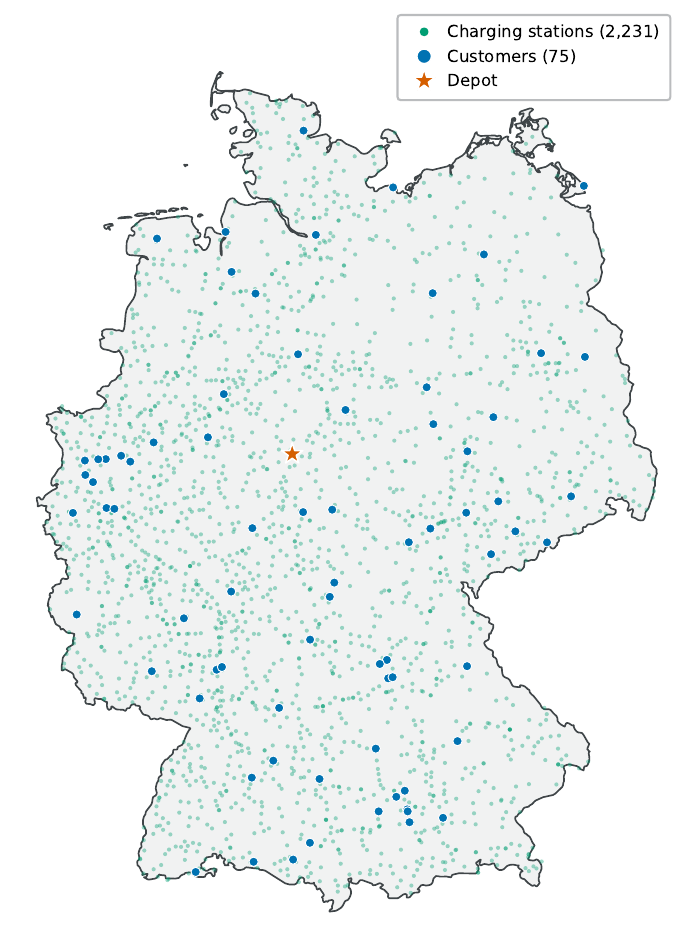}
    \caption{Example problem instance used in the experimental evaluation.}
    \label{fig:sample_instance}
\end{figure}

Figure \ref{fig:sample_instance} shows the total number of charging stations without any instance-specific filtering, and the customers and depot of an example problem instance with 75 customers.

The Open Source Routing Machine (OSRM) \citep{luxen2011real} is used for computing the travelling times and distances between each pair of nodes. Let \(d_{ij}\) and \(v_{ij}\) denote the resulting distance and average speed for an edge, respectively, with \(v_{ij}\) clipped to the interval \([1,80]\) km/h. Energy consumption is affected by vehicle speed, among other operational and road characteristics \citep{basso2019energy}. 

Because the detailed speed profiles, road gradients, and acceleration patterns required by a physical consumption model are not available for every edge, we employ a stylized speed-dependent energy-consumption function:
\[
e_{ij}=d_{ij}\left(0.8+0.00042(v_{ij}-70)^2\right),
\]
where \(e_{ij}\) is measured in kWh when \(d_{ij}\) is expressed in km. This model assumes a minimum consumption of \(0.8\) kWh/km at 70 km/h and a quadratic increase as the average speed moves away from this value.

We consider a fleet of battery-electric heavy-duty vehicles, each with a battery capacity of 850~kWh, a minimum admissible state of charge of \(y_{\min} = 50~\text{kWh}\), and a maximum charging power of 1,200~kW. These values are motivated by the specifications of the Tesla Semi, for which a range of 500 miles and an energy consumption of approximately 1.7~kWh/mile imply a battery capacity of approximately 850~kWh, while charging at up to 1.2~MW is supported \citep{teslaSemi2026}. For simplicity, we assume that this maximum charging power is available at every charging station. 

The fleet is homogeneous and effectively unlimited. Customer demands and vehicle capacity are expressed in standardized units, and the capacity of each vehicle is set to 100. Consequently, the demand values introduced above should be interpreted relative to this standardized capacity rather than as physical units of cargo.

The operational-cost coefficients are set to \(C_e=0.4\)~\(\text{€}/\text{kWh}\), \(C_{\mathrm{time}}=8\)~\(\text{€}/\text{h}\), and \(C_{\mathrm{veh}}=200\)~€ per route. These are experimental baseline weights used to define the trade-offs between the components of the objective and should not be interpreted as calibrated estimates of commercial operating costs. Unless stated otherwise, charging-related emissions are assigned a weight of \(p_{\mathrm{CO2}}=10\)~\(\text{€}/\mathrm{kgCO_2eq}\). Because this value is substantially higher than conventional market carbon prices, \(p_{\mathrm{CO2}}\) should be interpreted as an experimental penalty coefficient controlling the trade-off between operational cost and charging-related emissions, rather than as a literal market price for carbon.

\subsection{Data}

Problem instances are assigned routing start dates between 2023 and 2025. To predict grid carbon intensity over each 72-hour routing horizon, we combine electricity grid observations with weather forecasts at hourly resolution. The electricity grid features comprise the life-cycle carbon intensity and total load of the German grid, while the weather features describe forecast solar irradiance and wind conditions.

The carbon intensity data were obtained from the Electricity Maps Data Portal \citep{electricitymaps2026}, while load observations were obtained directly from the ENTSO-E Transparency Platform \citep{entsoe_transparency}. Electricity Maps derives consumption-based carbon intensities from electricity-generation and cross-border-flow data, whose underlying source for Germany is also ENTSO-E. For each routing instance, carbon-intensity and load observations from the 24 hours preceding the routing start are used as lagged predictive features, as described in the following subsection.

Weather forecasts were obtained from the \(0.25^\circ\times0.25^\circ\) Global Forecast System (GFS) grids produced by the National Centers for Environmental Prediction \citep{ncar_gdex_dataset_d084001}. Specifically, we extracted surface downward short-wave radiation and the horizontal wind components at 100 m from the pgrb2.0p25 product using the Herbie software package \citep{blaylock_herbie}. The GFS forecast horizon covers the three days considered by the routing instances.

To construct Germany-level weather features, the forecast associated with each grid cell is weighted according to the renewable-generation capacity installed within that cell. Solar irradiance is weighted by installed solar capacity, while the two wind components are weighted by installed wind capacity. Installation locations and capacities were obtained from the German Market Master Data Register (MaStR) \citep{bnetza_mastr} and assigned to the corresponding \(0.25^\circ\times0.25^\circ\) grid cells. Because the MaStR capacity data are available as quarterly snapshots, each forecast date is associated with the temporally closest available snapshot.

The resulting dataset covers 2023--2025 at hourly resolution and contains the grid's life-cycle carbon intensity and load, the capacity-weighted solar-irradiance forecasts, and the capacity-weighted forecasts of the two horizontal wind components at 100~m. The magnitude of the capacity-weighted wind vector is additionally derived from its two components and also included.

\subsection{Carbon Intensity Model}

As introduced in Section~\ref{subsection:training_loop} when presenting the training loop, we use a linear model to predict the grid carbon intensity over the routing horizon. This choice keeps the predictive component transparent and enables a controlled comparison between \gls{pfl} and \gls{dfl}. Moreover, it permits the direct computation of the subgradient given in Equation~(\ref{eq:subgradient}). Investigating whether richer predictive models yield further downstream improvements is left for future work.

We forecast the carbon intensity of the grid at hourly resolution and accordingly use hourly time buckets in the routing objective. Let \(H=\{0,\ldots,71\}\) denote the hourly forecast periods and let \(t_h=t_s+h\), where \(t_s\) is the routing start time. The model produces 72 predictions covering
the interval \([t_s,t_s+72\text{ h})\).

For each time step \(h\in H\), we define
\[
    \tau_h=t_s-24\text{ h}+(h\bmod 24) \ \text{h}.
\]
Thus, \(\tau_h\) is the corresponding hour from the day before the routing start time. The feature vector for period \(h\) contains:
\begin{enumerate}
    \item the observed grid carbon intensity at time \(\tau_h\).
    \item the observed grid load at time \(\tau_h\).
    \item the capacity-weighted solar irradiance forecast at time \(t_h\).
    \item the capacity-weighted forecasts for both wind components and wind-speed magnitude at 100 m at time \(t_h\). 
    \item cyclical encodings of the hour of the day, day of the week, and day of the year associated with \(t_h\).
\end{enumerate}

For a temporal variable \(q\) with period \(P\), its cyclical encoding is
defined as
\[
    \operatorname{cyc}(q;P)
    =
    \left(
        \sin\left(\frac{2\pi q}{P}\right),
        \cos\left(\frac{2\pi q}{P}\right)
    \right).
\]
As is standard in time-series prediction, we use periods \(24\), \(7\), and \(365\) for the hour of the day, day of
the week, and day of the year, respectively. This results in six features encoding time.

Before fitting the model,
the six non-cyclical features are standardized using means and standard
deviations computed from the training period. The cyclical features are
left unchanged.

Let \(z_h\in\mathbb{R}^{12}\) denote the resulting feature vector. Including an intercept gives \(x_h=(z_h^\top,1)^\top\in\mathbb{R}^{13}\). Stacking the feature vectors for all \(d=72\) periods produces a matrix \(X\in\mathbb{R}^{d\times13}\). The predicted carbon-intensity vector is then
\[
    \hat{c}=X\theta,
    \qquad
    \theta\in\mathbb{R}^{13},
    \qquad
    \hat{c}\in\mathbb{R}^{d}.
\]

The total rolling-horizon dataset contains 72 hourly observations for each available forecast origin, with one forecast origin per day at 00:00 UTC. Before removing missing values, the period 2023--2025 contains \((365+366+365)\cdot72=78{,}912\) forecast-origin and lead-time combinations. After removing missing observations, 78,765 remain: 52,560 from 2023--2024, which are used for fitting the \gls{pfl} model, and 26,205 from 2025, which are reserved for predictive evaluation. The \gls{dfl} training and validation instances are also drawn exclusively from 2023--2024, while downstream evaluation uses problem instances from 2025.

When the \gls{pfl} or \gls{dfl}-trained linear model produces a negative carbon-intensity prediction, that prediction is clipped to zero before it is supplied to the \gls{lns}. The same clipping is applied when computing the predictive-performance metrics reported below.

\subsection{Results}

The goal of the experimental evaluation is twofold. First, we compare the proposed \gls{dfl} pipeline against standard \gls{pfl}. In the latter, we train the linear regression in a supervised learning manner using Ordinary Least Squares (OLS), with the real carbon intensity values as training labels. Second, we assess whether incorporating the carbon intensity term in the objective function leads to better solutions, as indicated by lower emissions.

\subsubsection{Training and Model Selection}

The \gls{dfl} weights are initialized using the coefficients of the \gls{pfl} linear model fitted by OLS on the 2023--2024 hourly observations. The decision-focused training then iteratively updates these parameters using the normalized approximate \spoplus{} directions described in Section~\ref{subsection:training_loop}.

\Gls{dfl} training uses 96 problem instances for training and 24 for validation. For the training set, 48 instances are sampled uniformly at random from each year, while the validation set contains 12 randomly sampled instances per year. Routing start dates are sampled uniformly over the valid days of the corresponding year, subject to the availability of a complete 72-hour forecast horizon. Each instance contains between 50 and 100 customers, with the number of customers sampled uniformly at random. The train-validation split is hence not stratified by routing-start month or number of customers.

The batch size is set to 2, while the learning rate is set to \(\eta=25\). However, this value is not the effective per-instance step size because every instance subgradient is divided by \(d\alpha_i\), where \(d=72\) and \(\alpha_i\) is the total energy charged by its reference solution. Across the training instances, the resulting factor \(\eta/(d\alpha_i)\) ranged from \(1.36\times10^{-5}\) to \(4.76\times10^{-5}\), with a median of \(2.25\times10^{-5}\).

Regarding the \gls{lns} executions during training, we use a time limit of 40 seconds for training executions and 60 seconds for validation ones, which aims to match the 60-second time limit used during evaluation. Moreover, we set the number of customers removed by the Shaw destroy operator, \(q\), to 10, the determinism parameter \(p\) to 4, and all relatedness weights $\phi_1, \phi_2, \phi_3$ to one, giving the same weight to each standardized relatedness measure. Finally, we set its $k$-repair operator's $k$ parameter to 2, and the initial temperature \(T_\text{init}\) and cooling factor \(c_f\) of its simulated annealing to 100 and 0.97, respectively.

To produce the reference solutions obtained under the true carbon intensity used in the training loop, we execute the \gls{lns} once with a time limit of three hours using the same parameters as above, except for the simulated-annealing parameters, which are adjusted to the long time limit, by setting $T_\text{init} = 1000$ and $c_f = 0.99995$.

\begin{figure}[pos=htbp]
    \centering
    \includegraphics[width=0.8\textwidth]{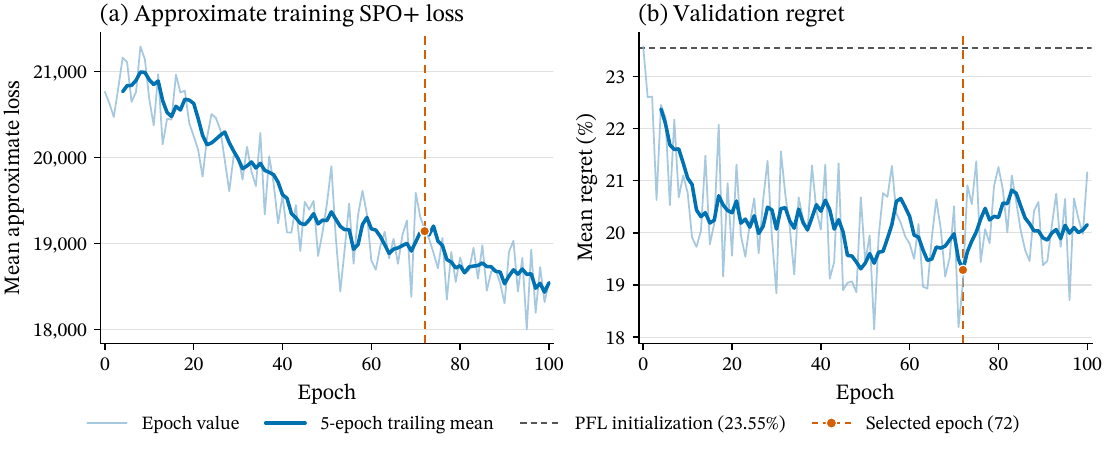}
    \caption{Training dynamics over epochs. The figure on the left presents the evolution of the approximate \spoplus{} loss, while the figure on the right presents the evolution of the validation regret percentage.}
    \label{fig:training_dynamics}
\end{figure}

Training was run for 100 epochs and required approximately 64 hours in total, with each epoch taking about 37 minutes on average. Overfitting was evaluated on the validation set after every epoch, and the five-epoch trailing average regret was used to decide which weights to retain. This was done with the goal of smoothing the epoch-by-epoch variability produced by the stochasticity of the \gls{lns}. Following this, the weights from epoch 72 were the ones retained as they led to the lowest five-epoch trailing average validation regret.

The clipping applied in line 6 of Algorithm~\ref{alg:Training_loop}, aimed at avoiding negative carbon intensity predictions being fed into the \gls{lns}, was infrequent during the selected training run. Across all constructed transformed carbon-intensity vectors, merely \(0.0047\%\) of their components were negative and therefore clipped to zero. At least one component was clipped in \(0.084\%\) of the \spoplus{} \gls{lns} solves. The difference between these percentages indicates that, when clipping occurred, it generally affected only a small part of the 72-period transformed vector. Thus, the projection rarely modified the input supplied to the downstream solver, although it still introduced an additional approximation whenever active.

Figure~\ref{fig:training_dynamics} shows the approximate \spoplus{} training loss and validation regret over the training epochs. The reported training loss would coincide with the theoretical \spoplus{} loss if both solutions involved were obtained from exact optimization oracles and the transformed cost vector was not clipped. To avoid additional \gls{lns} executions, each epoch value averages losses evaluated at the parameter values current during its successive mini-batches, rather than at a single parameter vector at the end of the epoch. The reported epoch average should therefore be interpreted as an online diagnostic rather than as the empirical loss at a single parameter vector. Validation regret is the average percentage excess cost under the true carbon intensity of the validation solutions relative to the precomputed reference solutions.

The initial validation regret, computed with the warm-started weights from \gls{pfl}, is represented by a gray horizontal dashed line. As can be observed, most learning happens within the first 20 epochs, as it is when validation regret decreases the most. We also plot the 5-epoch trailing mean of both measures, which in the case of the validation value was used for selecting the trained weights.

All training and evaluation runs were executed on Intel Xeon Gold 6338 ("Ice Lake") processors with 32 cores per CPU (64 cores total per node) running at 2.00 GHz. In particular, the two \gls{lns} executions required for processing each batch were done in parallel, while parallelization with 8 concurrent solves was used for the 24 \gls{lns} validation executions done at the end of every epoch. For solving the postprocessing model, presented in Section \ref{subsection:postprocessing} and applied at the end of every \gls{lns} execution, version 13.0 of the commercial solver Gurobi \citep{gurobi} was employed, with 8 threads reserved for every solve.

\subsubsection{Predictive Performance}

Before evaluating downstream decision quality, we compare the pointwise
predictive performance of the \gls{pfl} and \gls{dfl}-trained models. As a baseline, we also
consider a naive 24-hour persistence forecast that repeats the carbon
intensity observed during the 24 hours preceding the routing start over
the three-day routing horizon. The evaluation uses all available
forecast-origin and lead-time combinations from 2025, resulting in
26,205 observations.

Table~\ref{tab:predictive_performance} reports the mean absolute error
(MAE) and root mean squared error (RMSE) of the three approaches. Both metrics are expressed in grams of carbon-dioxide equivalent per kilowatt-hour ($\mathrm{gCO_2eq/kWh}$).

\begin{table}[pos=htbp]
    \centering
    \caption{Predictive performance on the 2025 evaluation data. Lower
    values indicate better pointwise predictive accuracy.}
    \label{tab:predictive_performance}
    \begin{tabular}{lrr}
        \hline
        Model & MAE & RMSE \\
        \hline
        \gls{pfl}                         & 57.70  & 72.46  \\
        \gls{dfl}                         & 402.70 & 410.71 \\
        Naive 24-hour persistence   & 90.88  & 120.03 \\
        \hline
    \end{tabular}
\end{table}

The \gls{pfl} model achieves the best pointwise predictive performance,
reducing MAE by 36.5\% and RMSE by 39.6\% relative to the naive
persistence baseline. This indicates that the linear model extracts
useful predictive information from the lagged grid variables, weather
forecasts, and cyclical time features beyond merely repeating the
carbon-intensity profile of the preceding day.

The \gls{dfl} model produces substantially larger pointwise errors. This is
not unexpected, since its parameters are optimized using the \spoplus{} loss
rather than a conventional prediction loss. Consequently, \gls{dfl} does not
directly penalize deviations between predicted and observed carbon
intensities. Its outputs are therefore better interpreted as
decision-oriented cost coefficients supplied to the routing algorithm
than as calibrated carbon-intensity forecasts. Whether sacrificing
predictive accuracy is beneficial must consequently be assessed through
the quality of the resulting routing decisions, as considered in the
following subsection.

\begin{figure}[pos=t]
    \centering
    \includegraphics[width=0.8\textwidth]{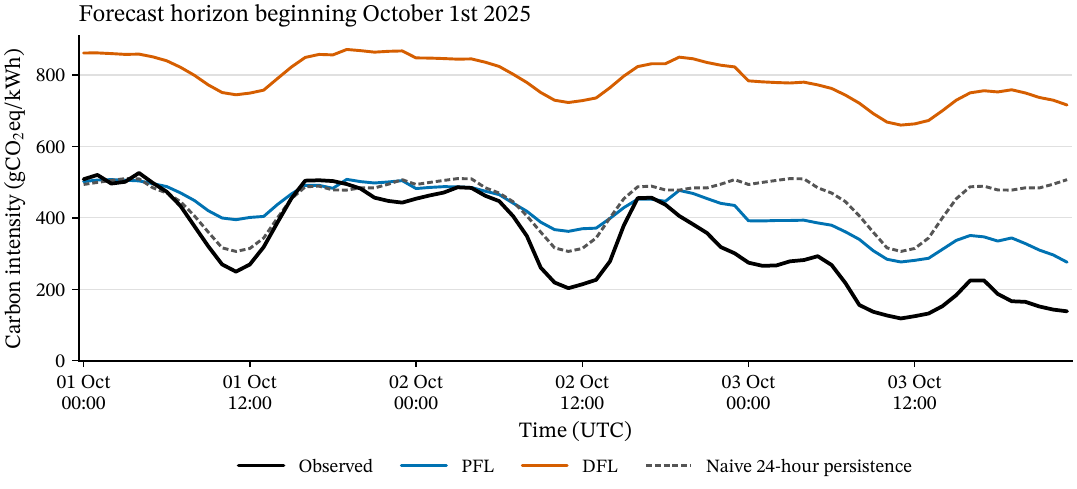}
    \caption{Observed carbon intensity and predictions produced by the \gls{pfl},
    \gls{dfl}, and naive 24-hour persistence models over the 72-hour routing
    horizon beginning on 1 October 2025.}
    \label{fig:example_predictions}
\end{figure}

Figure~\ref{fig:example_predictions} illustrates the predictions for
the 72-hour routing horizon beginning on October 1, 2025. In this
example, the \gls{pfl} prediction follows the observed carbon intensity more
closely than the naive baseline. The \gls{dfl} prediction exhibits temporal
variation but is shifted substantially upward, consistent with its
large aggregate MAE and RMSE.

\subsubsection{Downstream Decision Quality}

We evaluate downstream decision quality using 200 problem instances with routing start dates sampled uniformly over 2025. Each \gls{lns} execution is assigned a 60-second time limit and uses the baseline parameter configuration employed during training. The only parameter modified is the simulated-annealing cooling factor, which is increased from 0.97 to 0.98 because the evaluation time limit is longer than the 40-second limit used for the \gls{lns} executions during training.

We compare the \gls{dfl}-trained predictor with the \gls{pfl}-trained predictor, the naive 24-hour persistence forecast introduced above, and a perfect-information heuristic benchmark that supplies the true carbon intensity to the \gls{lns}. To account for the stochasticity of the \gls{lns}, each method is executed ten times on every instance, resulting in 2,000 executions per method. Regardless of the carbon-intensity forecast supplied during optimization, every resulting solution is subsequently evaluated using the true carbon intensity. Consequently, the comparison measures realized downstream cost under a common set of observed carbon-intensity values.

The 60-second limit applies only to the \gls{lns}. The exact fixed-route postprocessing is applied once to the final solution after this limit and requires an additional 0.287 seconds per solution on average. The corresponding median and 95th percentile are 0.247 and 0.58 seconds, respectively.

For each method and problem instance, we first average the realized total cost over its ten executions. For baseline \(b\), the instance-level improvement of \gls{dfl} is
\[
I_b(i)=100\frac{\bar C_b(i)-\bar C_{\mathrm{DFL}}(i)}
{\bar C_b(i)}.
\]
Here, \(\bar C_b(i)\) and \(\bar C_{\mathrm{DFL}}(i)\) are the mean costs for instance \(i\) using baseline \(b\) and using the \gls{dfl}-trained predictor, respectively. Hence, positive values indicate that the \gls{dfl}-trained predictor achieves a lower mean realized cost than baseline \(b\). Statistical summaries therefore treat the 200 instance-level mean improvements, rather than the 2,000 individual executions per method, as the observations.

\begin{figure}[pos=t]
    \centering
    \includegraphics[width=0.8\textwidth]{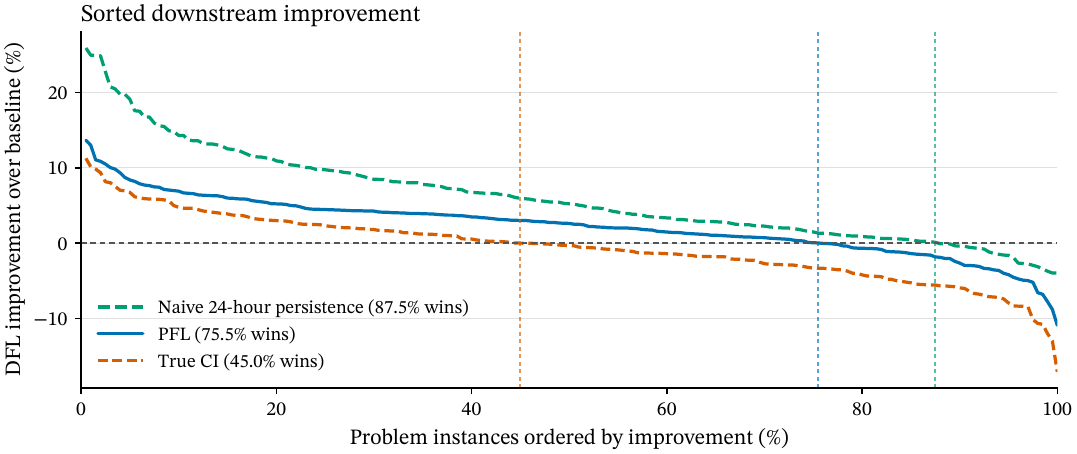}
    \caption{Sorted downstream improvement of the \gls{dfl}-trained predictor over the \gls{pfl}-trained predictor, the naive 24-hour persistence forecast, and the true-carbon-intensity benchmark.}
    \label{fig:improvement_curve}
\end{figure}

Figure~\ref{fig:improvement_curve} presents the sorted instance-level improvement of \gls{dfl} relative to each benchmark. The instances are ordered independently for each comparison; consequently, the same horizontal position does not generally correspond to the same problem instance across the three curves. Positive values indicate that \gls{dfl} achieves a lower mean realized cost than the corresponding benchmark. The vertical dashed lines mark the fraction of instances for which this improvement is positive.

Relative to \gls{pfl}, \gls{dfl} reduces realized cost by 2.30\% on average, with a median improvement of 2.56\%, and achieves a lower average cost on 75.5\% of the instances. Its advantage is larger relative to the naive persistence forecast, with a mean improvement of 6.16\% and a win rate of 87.5\%. By contrast, \gls{dfl} is 0.59\% more costly on average than the perfect-information heuristic benchmark and outperforms it on 45\% of the instances. This benchmark has access to the true carbon intensity but still relies on the same stochastic, time-limited \gls{lns} and is therefore not an optimal lower bound. Together with Table \ref{tab:predictive_performance}, these results show that the lower pointwise accuracy of \gls{dfl} nevertheless translates into better downstream decisions than the other forecast-based approaches.

\begin{figure}[pos=htbp]
    \centering
    \includegraphics[width=\textwidth]{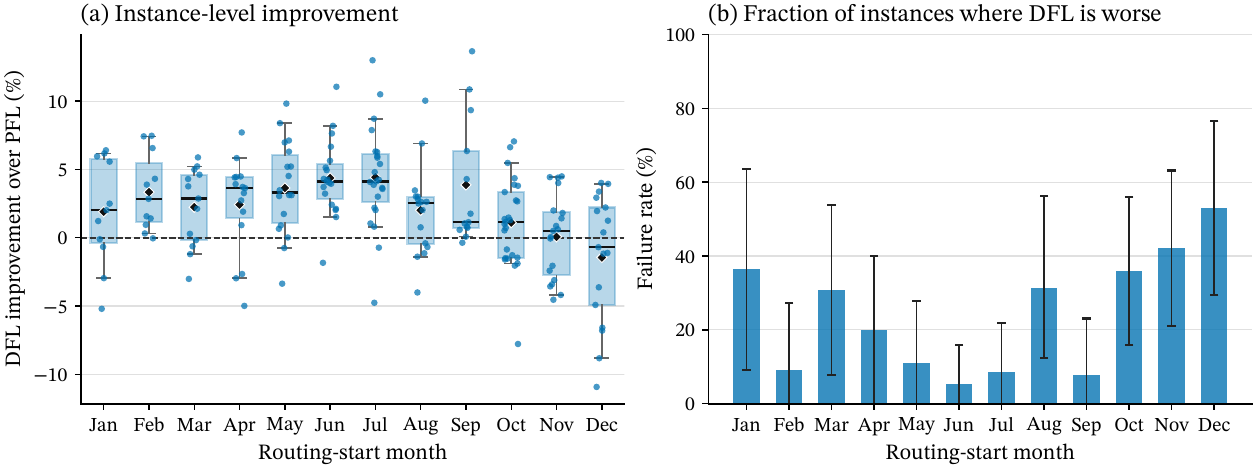}
    \caption{Downstream solution improvement of the \gls{dfl}-trained linear regression over the \gls{pfl}-trained one, with the results organized by the month of the instances' start date.}
    \label{fig:monthly_improvement}
\end{figure}

Figure~\ref{fig:monthly_improvement} examines how the relative performance of \gls{dfl} and \gls{pfl} varies with the routing-start month. Figure~\ref{fig:monthly_improvement}(a) shows the distribution of instance-level improvements within each month. Each point represents the ten-run mean improvement for one problem instance; center lines indicate medians, diamonds indicate means, boxes show interquartile ranges, and whiskers extend from the 5th to the 95th percentile. Figure~\ref{fig:monthly_improvement}(b) reports the corresponding failure rate, defined as the proportion of instances for which \gls{dfl} produces a higher mean realized cost than \gls{pfl}. The error bars indicate 95\% percentile bootstrap confidence intervals, computed from 20,000 resamples with replacement within each month.

Performance varies with the routing month. Mean \gls{dfl} improvement remains positive from January through October but becomes approximately zero in November and negative in December. Across November--December, the mean improvement is \(-0.69\%\), compared with \(2.93\%\) over January--October, while the corresponding failure rates are 47.5\% and 19.6\%. Because this seasonal analysis was conducted post hoc and the number of test instances differs across months, it should be interpreted as exploratory evidence rather than as a statistical conclusion.

\subsubsection{Small-scale Exact Validation}

To assess the solution quality of the \gls{lns} against proven-optimal solutions, we generated a deliberately small benchmark consisting of ten instances with each of three, four, and five customers. All locations are in the region of Saxony, and each instance has a 24-hour routing horizon. Depot and
customer locations were sampled from the same preprocessed OpenStreetMap locations used for the main instances. For each instance, we retained charging stations within 0.15~km of at least one
OSRM route connecting a pair of service locations, comprising the depot and customers, and subsequently removed stations without a feasible incoming connection. This resulted in averages of 5.6, 7.5, and 9.3 charging stations for the three, four, and five customer instances, respectively.

The standardized vehicle load capacity remained \(Q_{\max}=100\). To preserve the relevance of the capacity constraint despite the small number of customers, aggregate demand was sampled uniformly between 140 and 180 and distributed evenly across customers. Regarding time windows, each customer was given a \(2/3\) probability of being assigned a flexible time window covering the whole 24 hours, and \(1/3\) probability of being assigned a working-hours time window, with the opening and closing times sampled uniformly at random from hours 5--9 and 17--21, respectively. The battery capacity was reduced from 850 to 300~kWh so that intermediate charging remained relevant. Both methods were supplied with the true carbon-intensity trajectory. Gurobi was given a ten-hour time limit and used two copies of each physical charging station, while the \gls{lns} was executed ten times per instance with a 60-second limit. Consequently, the optimality statements below refer to the \gls{milp} containing two copies of each physical station.

\begin{table}[pos=htbp]
    \centering
    \caption{Small-scale comparison of the \gls{lns} with the exact \gls{milp}.}
    \label{tab:exact_benchmark}
    \begin{tabular}{rrrrr}
        \hline
        Customers & Mean stations & Optimality Proven & Mean gap (\%) & Median gap (\%) \\
3 & 5.6 & 10/10 & 1.46 & 0.77 \\
4 & 7.5 & 5/10  & 2.48 & 1.43 \\
5 & 9.3 & 0/10  & --   & --   \\
        \hline
    \end{tabular}
\end{table}

Table~\ref{tab:exact_benchmark} presents the resulting comparison, with each instance-level gap calculated using the mean objective over the ten \gls{lns} executions. Gurobi proved optimality for all three-customer instances, against which the \gls{lns} obtained a mean gap of 1.46\% and a median gap of 0.77\%. For the four-customer instances, it proved optimality for five; over this solved subset, the corresponding gaps were 2.48\% and 1.43\%. Finally, it could not prove optimality for any five-customer instance within the ten-hour limit. Although these instances are necessarily small, the results provide a direct validation of \gls{lns} solution quality where exact comparison is possible and illustrate the rapid computational growth of the \gls{milp}.

\subsubsection{Robustness to the \gls{lns} Configuration}

The main evaluation uses the same \gls{lns} configuration as \gls{dfl} training, except for the cooling factor, which is adjusted to account for the longer execution time. Furthermore, its 60-second time limit is identical to that used when evaluating validation regret during training. The learned predictor could therefore be specialized to the solver configuration used throughout training and model selection. To assess whether its advantage generalizes beyond this configuration, we reevaluate the fixed \gls{dfl} and \gls{pfl} predictors under alternative \gls{lns} settings. Notably, the predictors are not retrained in any of these experiments.

\begin{figure}[pos=t]
    \centering
    \includegraphics[width=0.8\textwidth]{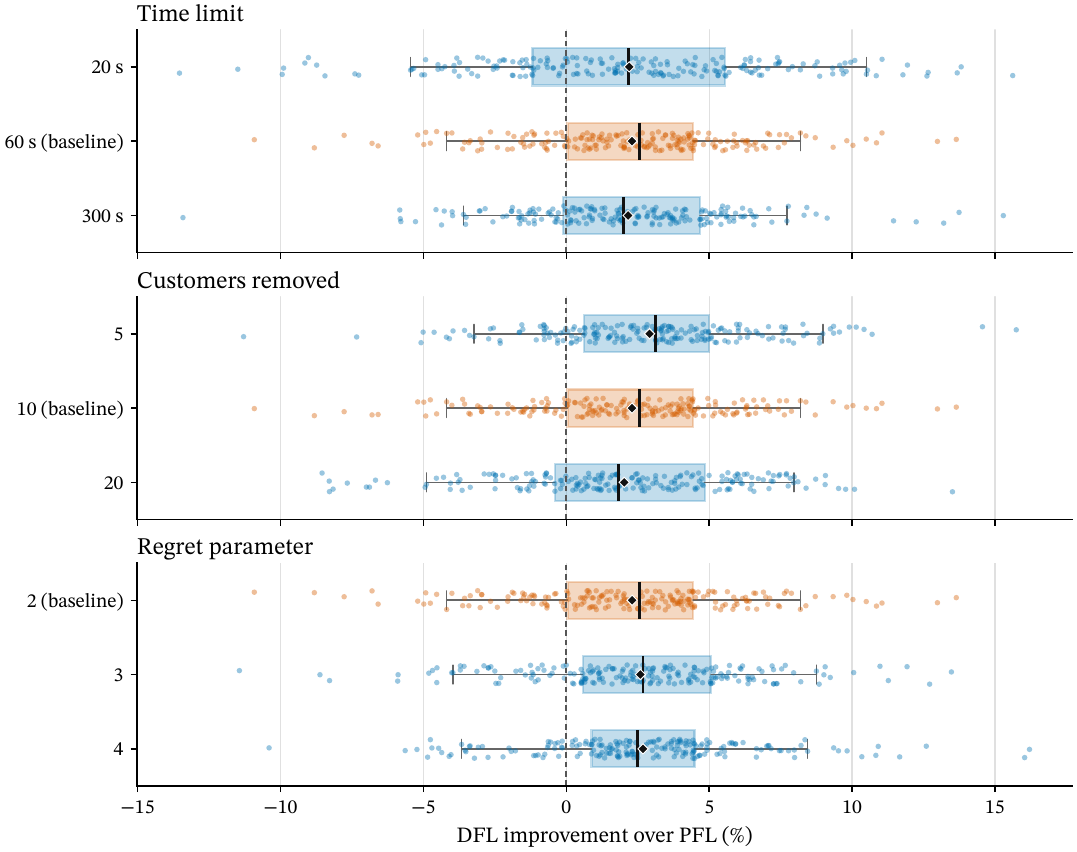}
    \caption{Downstream improvement of the \gls{dfl}-trained predictor over the \gls{pfl}-trained predictor under different \gls{lns} parameter configurations.}
    \label{fig:robustness_lns}
\end{figure}

Figure~\ref{fig:robustness_lns} presents the instance-level improvement of \gls{dfl} over \gls{pfl} under different evaluation configurations. We vary one of the \gls{lns} main parameters at a time while holding the remaining parameters at their baseline values. Specifically, we consider execution time limits of 20, 60, and 300 seconds; destroy sizes \(q\in\{5,10,20\}\); and regret parameters \(k\in\{2,3,4\}\). The baseline configuration uses a 60-second time limit, \(q=10\), and \(k=2\). When changing the time limit, the cooling factor is adjusted so that the simulated-annealing temperature follows a comparable cooling schedule over the duration of the search. Particularly, it is set to \(0.941192\) and \(0.995968\), for the 20 and 300 second time limits, respectively. Each point represents the ten-run mean improvement for one problem instance; center lines indicate medians, diamonds indicate means, boxes show interquartile ranges, and whiskers extend from the 5th to the 95th percentile.

The \gls{dfl} advantage remains positive across all tested \gls{lns} configurations. Mean improvement over \gls{pfl} ranges from 2.03\% to 2.91\% when varying the time limit, destroy size, or regret parameter. Particularly, increasing the time limit seems to reduce variability while the mean and median improvements remain positive. Increasing the number of customers removed seems to decrease the mean and median improvement. Finally, increasing the regret parameter does not result in noticeable changes in performance. These results therefore suggest that the learned predictor is not narrowly tied to the precise solver configuration used during training.

\subsubsection{Sensitivity to Carbon-cost Weight}
\begin{figure}[pos=t]
    \centering
    \includegraphics[width=0.8\textwidth]{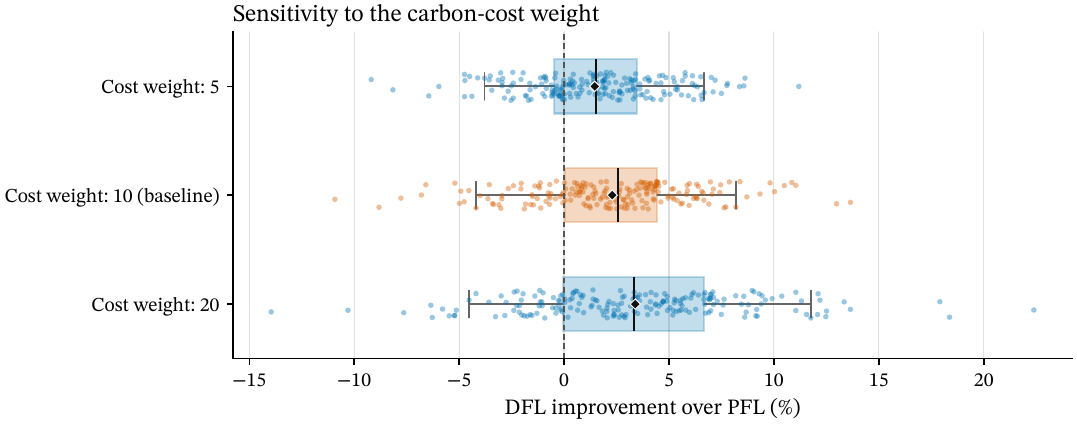}
    \caption{Downstream improvement of the \gls{dfl}-trained predictor over the \gls{pfl}-trained predictor for carbon-cost weights of 5, 10, and 20.}
    \label{fig:robustness_carbon}
\end{figure}

We also assess sensitivity to the carbon-cost weight \(p_{\mathrm{CO2}}\) used to value charging-related emissions. The training and evaluation use a carbon-cost weight of 10~\(\text{€}/\mathrm{kgCO_2eq}\), while the sensitivity analysis additionally considers carbon-cost weights of 5 and 20. As in the previous section, we evaluate these different configurations using the same trained predictor as in the evaluation, without retraining. Figure \ref{fig:robustness_carbon} shows a monotonic increase in the mean and median improvement of \gls{dfl} across the tested carbon-cost weights. Its mean improvement over \gls{pfl} is 1.52\%, 2.30\%, and 3.40\% for carbon-cost weights of 5, 10, and 20, respectively, while the corresponding fractions of instances on which \gls{dfl} performs better are 66\%, 75.5\%, and 74.5\%. As before, each point represents the ten-run mean improvement for one problem instance; center lines indicate medians, diamonds indicate means, boxes show interquartile ranges, and whiskers extend from the 5th to the 95th percentile.

\subsubsection{Effects of Considering the Carbon Cost}

Finally, we examine whether explicitly accounting for charging-related emissions changes the solutions produced by the \gls{lns} and postprocessing. We compare a carbon-aware configuration, which supplies the \gls{lns} with the \gls{dfl} predictor’s carbon-intensity forecasts and uses the baseline carbon-cost weight, with a carbon-blind configuration obtained by setting this weight to zero during optimization. All other \gls{lns} parameters and execution limits are identical to those used in the main evaluation. Similarly, both approaches are executed ten times on each of the 200 evaluation instances, and their results are first averaged within each instance.

Although the carbon-blind solutions ignore carbon intensity during optimization, they are evaluated afterward using the true carbon intensity and the baseline carbon-cost weight \(p_{\mathrm{CO_2}}=10\)~\(\text{€}/\mathrm{kgCO_2eq}\). This allows us to compare the charging-related emissions, non-carbon operating cost, and realized total cost of both approaches under a common objective. 

\begin{figure}[pos=t]
    \centering
    \includegraphics[width=0.8\textwidth]
        {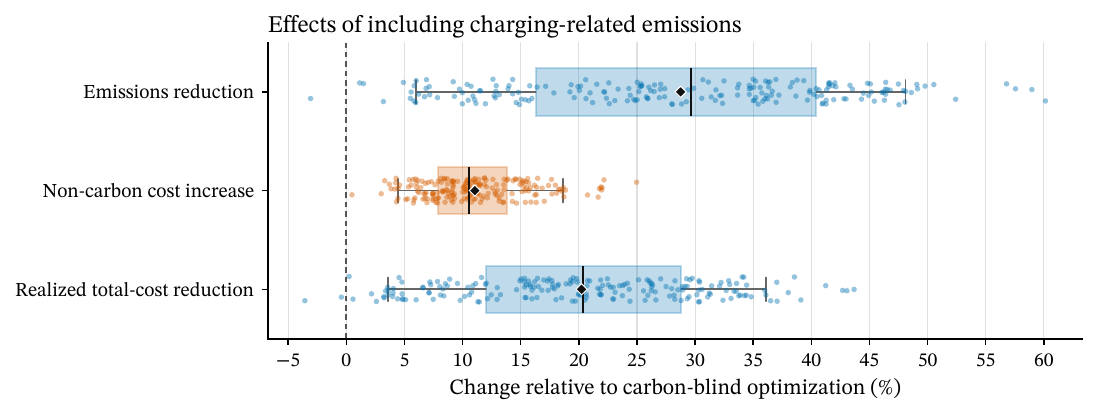}
    \caption{Effects of including charging-related emissions in the
    optimization objective.}
    \label{fig:no_carbon_effects}
\end{figure}

Figure~\ref{fig:no_carbon_effects} shows that explicitly considering charging-related emissions produces a substantial reduction in realized emissions. Relative to carbon-blind optimization, the carbon-aware solutions reduce emissions by 28.76\% on average across instances. These solutions produce lower emissions on 99.5\% of the evaluation instances. In absolute terms, average emissions decrease from 5,525.96 to 4,026.56~kg of \(\mathrm{CO_2eq}\) per instance.

As expected, this emissions reduction is accompanied by an increase in non-carbon operating cost. Specifically, mean non-carbon operating cost increases by 11.06\%, and it increases on every evaluation instance, illustrating the trade-off introduced by the carbon objective. Nevertheless, when both cost components are evaluated using the baseline carbon-cost weight, the carbon-aware solutions reduce realized total cost by 20.24\% on average, and a lower realized total cost is achieved on 99\% of the instances.

The carbon-aware solutions also require 2.05\% less charging energy on average. Because this reduction is considerably smaller than the 28.76\% reduction in emissions, the improvement cannot be explained by energy consumption alone. It also reflects changes in routing and charging timing that shift electricity consumption toward periods with lower realized carbon intensity. Overall, these results show that including charging-related emissions materially affects the resulting decisions and achieves substantial emission reductions, although at the cost of higher non-carbon expenditure.

The effect on fleet size is comparatively small. Under the baseline objective coefficients, the carbon-aware solutions use an average of 13.787 vehicles per instance, compared with 13.677 for the carbon-blind solutions, corresponding to an increase of 0.110 vehicles, or 0.8\%. These values should be interpreted conditionally on the experimental objective coefficients, since the relative weights assigned to vehicles, operating costs, and charging-related emissions determine the trade-off between fleet size and emission reduction. Nevertheless, the small increase in fleet size relative to the 28.76\% reduction in charging-related emissions indicates that the emission savings cannot be attributed primarily to the use of additional vehicles.

\FloatBarrier

\section{Conclusions and Future Work}\label{sec:conclusions}

This paper presented an \gls{evrp} formulation that incorporates charging-related emissions into its objective. The proposed formulation jointly optimizes routes, charging, and waiting decisions while accounting for temporal variation in grid carbon intensity. This is motivated by the substantial temporal variation in the electricity generation mix of some regions, since charging during periods of high grid carbon intensity can substantially increase charging-related emissions and, in sufficiently carbon-intensive grids, can even erode the lifecycle emissions advantage of electric vehicles over efficient hybrid alternatives.

We presented a decision-focused learning approach for tackling this problem. A linear predictor was trained using approximate \spoplus{} subgradients to predict the carbon intensity of the grid over the routing horizon, while a time-limited \gls{lns} with exact postprocessing was used as the downstream solver.

The experimental results demonstrate the usefulness of this carbon-aware formulation, as the solutions produced reduced emissions by \(28.76\%\) on average relative to carbon-blind optimization, while increasing the non-carbon cost by \(11.06\%\). Under the baseline objective coefficients, this trade-off resulted in a \(20.24\%\) reduction in realized total cost.

Moreover, the results also illustrate the distinction between predictive accuracy and downstream decision quality. Although the \gls{pfl}-trained model achieved substantially lower forecasting MAE and RMSE, the \gls{dfl}-trained predictor produced better solutions on 75.5\% of the test instances relative to \gls{pfl}. Its mean advantage also remained positive across the evaluated \gls{lns} configurations and carbon-cost weights, indicating robustness to the solver and carbon cost configurations. The exploratory monthly analysis nevertheless revealed weaker performance toward the end of the year, indicating that the benefits of decision-focused training were heterogeneous across the evaluation period.

The principal limitation is computational cost: the selected epoch (72) was reached after approximately 44 hours, while the complete 100-epoch training run required approximately 64 hours. To alleviate this issue, future work could consider further reducing the time limit assigned to the \gls{lns} during training and exploring the trade-off between this time limit and downstream performance. Another option would be to construct a surrogate model to reduce the number of \gls{lns} and postprocessing executions.

Future work could also investigate stratified or seasonally balanced training samples and examine whether the observed seasonal deterioration is caused by limited training coverage, changes in grid behavior, or characteristics of the routing instances. Also, richer predictive models could be considered. Finally, further extensions could consider routes spanning several electricity bidding regions, with region-specific carbon-intensity predictions.  

\section*{Acknowledgements}

The computations were enabled by resources provided by Chalmers
e-Commons at Chalmers. The computational experiments were performed
on the Vera computing cluster under project C3SE2026-1-13.

\printcredits

\bibliographystyle{cas-model2-names}

\bibliography{bib}

@article{elmachtoub2022smart,
  title   = {{Smart ``Predict, then Optimize''}},
  author  = {Elmachtoub, Adam N. and Grigas, Paul},
  journal = {Management Science},
  volume  = {68},
  number  = {1},
  pages   = {9--26},
  year    = {2022},
  doi     = {10.1287/mnsc.2020.3922}
}

@article{mandi2024decision,
  title   = {{Decision-Focused Learning: Foundations, State of the Art, Benchmark and Future Opportunities}},
  author  = {Mandi, Jayanta and Kotary, James and Berden, Senne and Mulamba, Maxime and Bucarey, Victor and Guns, Tias and Fioretto, Ferdinando},
  journal = {Journal of Artificial Intelligence Research},
  volume  = {80},
  pages   = {1623--1701},
  year    = {2024},
  doi     = {10.1613/jair.1.15320}
}

@article{singh2024ensuring,
  title   = {{Ensuring greenhouse gas reductions from electric vehicles compared to hybrid gasoline vehicles requires a cleaner U.S. electricity grid}},
  author  = {Singh, Madalsa and Yuksel, Tugce and Michalek, Jeremy J. and Azevedo, In{\^e}s M. L.},
  journal = {Scientific Reports},
  volume  = {14},
  number  = {1},
  pages   = {1639},
  year    = {2024},
  doi     = {10.1038/s41598-024-51697-1}
}

@article{su2025optimizing,
  title   = {{Optimizing carbon footprint in long-haul heavy-duty E-Truck transportation}},
  author  = {Su, Junyan and Lin, Qiulin and Chen, Minghua},
  journal = {Nature Communications},
  volume  = {16},
  number  = {1},
  pages   = {9562},
  year    = {2025},
  doi     = {10.1038/s41467-025-64792-2}
}

@article{kucukoglu2021electric,
  title   = {{The electric vehicle routing problem and its variations: A literature review}},
  author  = {Kucukoglu, Ilker and Dewil, Reginald and Cattrysse, Dirk},
  journal = {Computers \& Industrial Engineering},
  volume  = {161},
  pages   = {107650},
  year    = {2021},
  doi     = {10.1016/j.cie.2021.107650}
}

@incollection{jarvis2022smart,
  title     = {{Smart and Sustainable Scheduling of Charging Events for Electric Buses}},
  author    = {Jarvis, Padraigh and Climent, Laura and Arbelaez, Alejandro},
  booktitle = {{EU} Cohesion Policy Implementation: Evaluation Challenges and Opportunities},
  series    = {Springer Proceedings in Political Science and International Relations},
  pages     = {121--129},
  year      = {2022},
  publisher = {Springer Nature Switzerland},
  doi       = {10.1007/978-3-031-18161-0_8}
}

@article{miyabe2025low,
  title   = {{Low-carbon routing and charging planning for electric freight trucks utilizing local surplus solar power}},
  author  = {Miyabe, Ryoji and Fujimoto, Yu and Hayashi, Yasuhiro},
  journal = {Journal of Energy Storage},
  volume  = {132},
  pages   = {117626},
  year    = {2025},
  doi     = {10.1016/j.est.2025.117626}
}

@techreport{europeancommission2026electrification,
  author      = {{European Commission}},
  title       = {{Communication from the Commission to the European Parliament, the Council, the European Economic and Social Committee and the Committee of the Regions: Electrification Action Plan}},
  institution = {European Commission, Directorate-General for Energy},
  type        = {COM},
  number      = {(2026) 595 final},
  year        = {2026},
  month       = jul,
  url         = {https://eur-lex.europa.eu/legal-content/EN/TXT/?uri=CELEX:52026DC0595}
}

@misc{eurostat2026transport,
  author       = {{Eurostat}},
  title        = {{Final energy consumption in transport - detailed statistics}},
  year         = {2026},
  howpublished = {Eurostat Statistics Explained},
  url          = {https://ec.europa.eu/eurostat/statistics-explained/index.php?title=Final_energy_consumption_in_transport_-_detailed_statistics},
  note         = {{Accessed: 2026-09-21}}
}

@incollection{pisinger2018large,
  title     = {{Large Neighborhood Search}},
  author    = {Pisinger, David and R{\o}pke, Stefan},
  booktitle = {Handbook of Metaheuristics},
  series    = {International Series in Operations Research \& Management Science},
  volume    = {272},
  pages     = {99--127},
  year      = {2018},
  publisher = {Springer International Publishing},
  doi       = {10.1007/978-3-319-91086-4_4}
}

@article{chu2023data,
  title   = {{Data-driven optimization for last-mile delivery}},
  author  = {Chu, Hongrui and Zhang, Wensi and Bai, Pengfei and Chen, Yahong},
  journal = {Complex \& Intelligent Systems},
  volume  = {9},
  number  = {3},
  pages   = {2271--2284},
  year    = {2023},
  doi     = {10.1007/s40747-021-00293-1}
}

@inproceedings{mandi2021data,
  title     = {{Data Driven {VRP}: A Neural Network Model to Learn Hidden Preferences for {VRP}}},
  author    = {Mandi, Jayanta and Canoy, Rocsildes and Bucarey, V{\'\i}ctor and Guns, Tias},
  booktitle = {27th International Conference on Principles and Practice of Constraint Programming ({CP} 2021)},
  series    = {Leibniz International Proceedings in Informatics ({LIPIcs})},
  volume    = {210},
  pages     = {42:1--42:17},
  year      = {2021},
  publisher = {Schloss Dagstuhl--Leibniz-Zentrum f{\"u}r Informatik},
  doi       = {10.4230/LIPIcs.CP.2021.42}
}

@article{tian2023prescriptive,
  title   = {{Prescriptive analytics for a maritime routing problem}},
  author  = {Tian, Xuecheng and Yan, Ran and Wang, Shuaian and Laporte, Gilbert},
  journal = {Ocean \& Coastal Management},
  volume  = {242},
  pages   = {106695},
  year    = {2023},
  doi     = {10.1016/j.ocecoaman.2023.106695}
}

@article{kirkpatrick1983optimization,
  title   = {{Optimization by Simulated Annealing}},
  author  = {Kirkpatrick, Scott and Gelatt, Jr., C. Daniel and Vecchi, Mario P.},
  journal = {Science},
  volume  = {220},
  number  = {4598},
  pages   = {671--680},
  year    = {1983},
  doi     = {10.1126/science.220.4598.671}
}

@inproceedings{amos2017optnet,
  title     = {{{OptNet}: Differentiable Optimization as a Layer in Neural Networks}},
  author    = {Amos, Brandon and Kolter, J. Zico},
  booktitle = {Proceedings of the 34th International Conference on Machine Learning},
  series    = {Proceedings of Machine Learning Research},
  volume    = {70},
  pages     = {136--145},
  year      = {2017},
  publisher = {PMLR},
  url       = {https://proceedings.mlr.press/v70/amos17a.html}
}

@inproceedings{domke2012generic,
  title     = {{Generic Methods for Optimization-Based Modeling}},
  author    = {Domke, Justin},
  booktitle = {Proceedings of the Fifteenth International Conference on Artificial Intelligence and Statistics},
  series    = {Proceedings of Machine Learning Research},
  volume    = {22},
  pages     = {318--326},
  year      = {2012},
  publisher = {PMLR},
  url       = {https://proceedings.mlr.press/v22/domke12.html}
}

@inproceedings{kotary2023backpropagation,
  title     = {{Backpropagation of Unrolled Solvers with Folded Optimization}},
  author    = {Kotary, James and Dinh, My H. and Fioretto, Ferdinando},
  booktitle = {Proceedings of the Thirty-Second International Joint Conference on Artificial Intelligence},
  pages     = {1963--1970},
  year      = {2023},
  doi       = {10.24963/ijcai.2023/218}
}

@inproceedings{wilder2019melding,
  title     = {{Melding the Data-Decisions Pipeline: Decision-Focused Learning for Combinatorial Optimization}},
  author    = {Wilder, Bryan and Dilkina, Bistra and Tambe, Milind},
  booktitle = {Proceedings of the {AAAI} Conference on Artificial Intelligence},
  volume    = {33},
  number    = {01},
  pages     = {1658--1665},
  year      = {2019},
  doi       = {10.1609/aaai.v33i01.33011658}
}

@inproceedings{ferber2020mipaal,
  title     = {{{MIPaaL}: Mixed Integer Program as a Layer}},
  author    = {Ferber, Aaron and Wilder, Bryan and Dilkina, Bistra and Tambe, Milind},
  booktitle = {Proceedings of the {AAAI} Conference on Artificial Intelligence},
  volume    = {34},
  number    = {02},
  pages     = {1504--1511},
  year      = {2020},
  doi       = {10.1609/aaai.v34i02.5509}
}

@inproceedings{domke2010implicit,
  title     = {{Implicit Differentiation by Perturbation}},
  author    = {Domke, Justin},
  booktitle = {Advances in Neural Information Processing Systems},
  volume    = {23},
  pages     = {523--531},
  year      = {2010},
  url       = {https://proceedings.neurips.cc/paper/2010/hash/6ecbdd6ec859d284dc13885a37ce8d81-Abstract.html}
}

@inproceedings{berthet2020learning,
  title     = {{Learning with Differentiable Pertubed Optimizers}},
  author    = {Berthet, Quentin and Blondel, Mathieu and Teboul, Olivier and Cuturi, Marco and Vert, Jean-Philippe and Bach, Francis},
  booktitle = {Advances in Neural Information Processing Systems},
  volume    = {33},
  pages     = {9508--9519},
  year      = {2020},
  url       = {https://proceedings.neurips.cc/paper/2020/hash/6bb56208f672af0dd65451f869fedfd9-Abstract.html}
}

@inproceedings{niepert2021implicit,
  title     = {{Implicit {MLE}: Backpropagating Through Discrete Exponential Family Distributions}},
  author    = {Niepert, Mathias and Minervini, Pasquale and Franceschi, Luca},
  booktitle = {Advances in Neural Information Processing Systems},
  volume    = {34},
  pages     = {14567--14579},
  year      = {2021},
  url       = {https://proceedings.neurips.cc/paper/2021/hash/7a430339c10c642c4b2251756fd1b484-Abstract.html}
}

@inproceedings{mulamba2021contrastive,
  title     = {{Contrastive Losses and Solution Caching for Predict-and-Optimize}},
  author    = {Mulamba, Maxime and Mandi, Jayanta and Diligenti, Michelangelo and Lombardi, Michele and Bucarey, Victor and Guns, Tias},
  booktitle = {Proceedings of the Thirtieth International Joint Conference on Artificial Intelligence},
  pages     = {2833--2840},
  year      = {2021},
  publisher = {International Joint Conferences on Artificial Intelligence Organization},
  doi       = {10.24963/ijcai.2021/390}
}

@inproceedings{mandi2022decision,
  title     = {{Decision-Focused Learning: Through the Lens of Learning to Rank}},
  author    = {Mandi, Jayanta and Bucarey, V{\'\i}ctor and Tchomba, Maxime Mulamba Ke and Guns, Tias},
  booktitle = {Proceedings of the 39th International Conference on Machine Learning},
  series    = {Proceedings of Machine Learning Research},
  volume    = {162},
  pages     = {14935--14947},
  year      = {2022},
  publisher = {PMLR},
  url       = {https://proceedings.mlr.press/v162/mandi22a.html}
}

@inproceedings{shah2022decision,
  title     = {{Decision-Focused Learning without Decision-Making: Learning Locally Optimized Decision Losses}},
  author    = {Shah, Sanket and Wang, Kai and Wilder, Bryan and Perrault, Andrew and Tambe, Milind},
  booktitle = {Advances in Neural Information Processing Systems},
  volume    = {35},
  pages     = {1320--1332},
  year      = {2022},
  publisher = {Curran Associates, Inc},
  doi       = {10.52202/068431-0097}
}

@article{hart1968formal,
  title   = {{A Formal Basis for the Heuristic Determination of Minimum Cost Paths}},
  author  = {Hart, Peter E. and Nilsson, Nils J. and Raphael, Bertram},
  journal = {{IEEE} Transactions on Systems Science and Cybernetics},
  volume  = {4},
  number  = {2},
  pages   = {100--107},
  year    = {1968},
  doi     = {10.1109/TSSC.1968.300136}
}

@manual{gurobi,
  author = {{Gurobi Optimization, LLC}},
  title  = {{Gurobi Optimizer Reference Manual}},
  year   = {2026},
  note   = {{Version 13.0}},
  url    = {https://docs.gurobi.com/projects/optimizer/en/current/}
}

@misc{openstreetmap,
  author       = {{OpenStreetMap contributors}},
  title        = {{OpenStreetMap}},
  year         = {2026},
  howpublished = {OpenStreetMap},
  url          = {https://www.openstreetmap.org/},
  note         = {{Accessed: 2026-09-21}}
}

@inproceedings{ester1996density,
  title     = {{A Density-Based Algorithm for Discovering Clusters in Large Spatial Databases with Noise}},
  author    = {Ester, Martin and Kriegel, Hans-Peter and Sander, J{\"o}rg and Xu, Xiaowei},
  booktitle = {Proceedings of the Second International Conference on Knowledge Discovery and Data Mining},
  pages     = {226--231},
  year      = {1996},
  publisher = {{AAAI} Press},
  url       = {https://www.aaai.org/Papers/KDD/1996/KDD96-037.pdf}
}

@misc{ncar_gdex_dataset_d084001,
  author    = {{National Centers for Environmental Prediction, National Weather Service, NOAA, U.S. Department of Commerce}},
  title     = {{{NCEP GFS} 0.25 Degree Global Forecast Grids Historical Archive}},
  publisher = {{NSF} National Center for Atmospheric Research},
  address   = {Boulder, CO},
  year      = {2015},
  doi       = {10.5065/D65D8PWK}
}

@misc{blaylock_herbie,
  author    = {Blaylock, Brian K.},
  title     = {{Herbie: Retrieve Numerical Weather Prediction Model Data (Version 2026.3.0)}},
  year      = {2026},
  publisher = {Zenodo},
  doi       = {10.5281/zenodo.18902673}
}

@misc{bnetza_mastr,
  author       = {{Bundesnetzagentur}},
  title        = {{Marktstammdatenregister ({MaStR}): Gesamtdatenexport}},
  year         = {{n.d.}},
  howpublished = {Bundesnetzagentur},
  url          = {https://www.marktstammdatenregister.de/MaStR/Datendownload},
  note         = {{Quarterly full-data XML snapshots}}
}

@misc{electricitymaps2026,
  author       = {{Electricity Maps}},
  title        = {{Historical carbon intensity data for Germany}},
  year         = {2026},
  howpublished = {Electricity Maps Data Portal},
  url          = {https://portal.electricitymaps.com/datasets/DE},
  note         = {15-minute life-cycle carbon-intensity data for 2023--2025; export version dated 2026-02-10}
}

@misc{entsoe_transparency,
  author       = {{European Network of Transmission System Operators for Electricity ({ENTSO-E})}},
  title        = {{{ENTSO-E} Transparency Platform}},
  year         = {{n.d.}},
  howpublished = {{ENTSO-E}},
  url          = {https://transparency.entsoe.eu/},
  note         = {{Accessed: 2026-02-10}}
}

@inproceedings{luxen2011real,
  title     = {{Real-time routing with {OpenStreetMap} data}},
  author    = {Luxen, Dennis and Vetter, Christian},
  booktitle = {Proceedings of the 19th {ACM SIGSPATIAL} International Conference on Advances in Geographic Information Systems},
  pages     = {513--516},
  year      = {2011},
  publisher = {{ACM}},
  doi       = {10.1145/2093973.2094062}
}

@misc{teslaSemi2026,
  author       = {{Tesla}},
  title        = {{Tesla Semi}},
  year         = {2026},
  howpublished = {Tesla},
  url          = {https://www.tesla.com/semi},
  note         = {{Accessed: 2026-08-03}}
}

@article{basso2019energy,
  title   = {{Energy consumption estimation integrated into the Electric Vehicle Routing Problem}},
  author  = {Basso, Rafael and Kulcs{\'a}r, Bal{\'a}zs and Egardt, Bo and Lindroth, Peter and Sanchez-Diaz, Ivan},
  journal = {Transportation Research Part D: Transport and Environment},
  volume  = {69},
  pages   = {141--167},
  year    = {2019},
  doi     = {10.1016/j.trd.2019.01.006}
}

\end{document}